\documentclass[12pt,letterpaper]{article}
\usepackage{amsmath,amssymb,amsthm,mathtools}
\usepackage[lf]{Baskervaldx}
\usepackage[bigdelims,vvarbb]{newtxmath}
\usepackage[cal=boondoxo]{mathalfa}

\usepackage{helvet}
\usepackage[headings]{fullpage}
\renewcommand{\thefootnote}{\fnsymbol{footnote}}
\numberwithin{equation}{section}
\usepackage{xcolor}
\usepackage[colorlinks=true,linkcolor=blue!65!black,citecolor=blue!65!black,urlcolor=blue!65!black]{hyperref}
\usepackage{microtype}
\newtheorem{theorem}{Theorem}[section]
\newtheorem{lemma}[theorem]{Lemma}
\newtheorem{proposition}[theorem]{Proposition}
\newtheorem{corollary}[theorem]{Corollary}
\theoremstyle{remark}
\newtheorem{remark}[theorem]{Remark}
\DeclareMathOperator{\diver}{div}
\DeclareMathOperator{\tr}{tr}
\DeclareMathOperator{\supp}{supp}
\newcommand{\R}{\mathbb R}
\makeatletter
\def\NativeDashInt#1{\mathchoice{\NativeDashIntAux\displaystyle\textstyle{#1}}{\NativeDashIntAux\textstyle\scriptstyle{#1}}{\NativeDashIntAux\scriptstyle\scriptscriptstyle{#1}}{\NativeDashIntAux\scriptscriptstyle\scriptscriptstyle{#1}}\!\int}
\def\NativeDashIntAux#1#2#3{{\setbox0=\hbox{$#1{#2#3}{\int}$}\vcenter{\hbox{$#2#3$}}\kern-.5\wd0}}
\providecommand{\fint}{\NativeDashInt-}
\makeatother
\newcommand{\avg}{\mathop{\fint}}
\newcommand{\norm}[1]{\left\lVert#1\right\rVert}

\newcommand{\dd}{\,dx}
\newcommand{\wm}{\mathcal W}
\newcommand{\AP}{\mathcal A_p}
\title{An effective improvement of weighted second-order regularity for scalar $p$-Laplace equations}
\author{Sunghoon Kim\thanks{Department of Mathematics, The Catholic University of Korea, 43 Jibong-ro, Bucheon-si, Gyeonggi-do, 14662, Republic of Korea. E-mail: \texttt{math.s.kim@catholic.ac.kr}.}}
\date{}
\hypersetup{pdfauthor={Sunghoon Kim},pdftitle={An effective improvement of weighted second-order regularity for scalar p-Laplace equations}}
\begin{document}
\maketitle
\markboth{Sunghoon Kim}{Effective weighted second-order regularity}
\begingroup
\renewcommand{\thefootnote}{}
\footnotetext{\emph{Mathematics Subject Classification (2020):} 35J92, 35B65, 35D30, 35J70.}
\endgroup
\begin{abstract}
We study the local second-order Sobolev regularity of scalar $p$-harmonic functions. For $n\geq3$ and $p>2$, we give an explicit positive gain $\sigma(n,p)$ below the power threshold $m=(p-1)(n-2)/(n-1)$: the field $|Du|^{(m-\sigma(n,p)-2)/2}Du$ belongs to $W^{1,2}_{\rm loc}$. We first establish a lower bound for the tangential derivative of the gradient magnitude on interior balls with a prescribed oscillation. Its quantitative form follows from a Riccati inequality for the directional derivative of the logarithmic gradient magnitude along unit normal trajectories. This estimate, a covering argument, and the equation satisfied by the gradient magnitude yield the weighted Hessian estimate. We also consider equations with bounded right-hand sides and establish a corresponding estimate above the natural gradient scale, together with a Sobolev estimate for a regularized flux. In the homogeneous case, we derive higher integrability and distributional identities across the critical set. In particular, the range of $p$ for which $u\in W^{2,2+\varepsilon_{\rm hi}}_{\rm loc}$ extends beyond $3+2/(n-2)$ by an explicitly specified positive amount in every dimension $n\geq3$. Assuming, in addition, the weighted estimate stated in an unrefereed manuscript of Jin, we obtain this regularity in dimension three for all finite $p>2$.
\par\medskip\noindent\textbf{Keywords.} $p$-Laplacian, weighted Hessian, nonlinear gradient, critical set, higher integrability.
\end{abstract}
\vskip 0.2truein
\setlength{\baselineskip}{17pt}

\section{Introduction}

Let $\Omega\subset\R^n$ be open, $n\geq3$, and $p>2$. We consider scalar weak solutions of
\begin{equation}\label{eq:Equation}
 -\diver\bigl(\AP(Du)\bigr)=f,\qquad
 \AP(\xi)=|\xi|^{p-2}\xi,\qquad f\in L^\infty_{\rm loc}(\Omega).
\end{equation}
Throughout, write
\begin{equation*}
 P=p-1,\qquad b=p-2,\qquad m=\frac{P(n-2)}{n-1},\qquad Q=|Du|.
\end{equation*}
For $s>0$ define the vector field $N_s=|Du|^{(s-2)/2}Du$, with value zero where $Du=0$. Its length is $Q^{s/2}$.

The weak differentiability of $N_s$ is closely related to the integrability of $Q^{s-2}|D^2u|^2$ on the set $\{Q>0\}$. For the homogeneous equation, Dong, Peng, Zhang and Zhou \cite[Theorem 1.1]{DFZZ} proved the local $W^{1,2}$ regularity of $|Du|^{(p-\gamma)/2}Du$ under the condition
\begin{equation*}
 \gamma<\min\left\{p+\frac{n}{n-1},\ 3+\frac{p-1}{n-1}\right\}.
\end{equation*}
For $n\geq3$ and $p>2$, the second number in the minimum is the smaller one. Setting $s=p-\gamma+2$, we obtain precisely $s>m$. Subsequently, Sarsa \cite[Theorems 1.1 and 3.1]{Sarsa} refined the underlying algebraic inequality and obtained, in particular, a local estimate for the derivative of $|Du|^{(p-2+a)/2}Du$ with an arbitrary constant vector subtracted from this field on the right-hand side. The restriction $a>-1-(p-1)/(n-1)$ is again $s>m$ when $s=p+a$. Thus for the dimensions and exponents considered here the two sufficient power ranges coincide.

Karppinen and Sarsa \cite[Theorem 1.1 and Corollary 1.2]{KS} extended these estimates and the corresponding higher integrability to the Orlicz setting. In their power-function specialization, the field is $|Du|^\beta Du$, and the sufficient condition $\beta>-1+(p-1)(n-2)/(2(n-1))$ becomes $s>m$ upon setting $s=2(\beta+1)$. The question of improving this range in higher dimensions is discussed in \cite[Section 3.3]{KS}. We will investigate whether the restriction $s>m$ can be relaxed. More precisely, we prove the following result without imposing an assumption on the critical set of $u$.

\begin{theorem}[Weighted second-order estimate]\label{thm:main}
Let $n\geq3$ and $p>2$. There exist $s_*=s_*(n,p)\in(0,m)$ and $C=C(n,p)$ such that, if $u\in W^{1,p}_{\rm loc}(\Omega)$ is $p$-harmonic, then $N_{s_*}\in W^{1,2}_{\rm loc}(\Omega)$. For concentric balls $B_r\subset B_R\Subset\Omega$,
\begin{equation}\label{eq:main-est}
 \int_{B_r\cap\{Q>0\}} Q^{s_*-2}|D^2u|^2\dd
 +\int_{B_r}|DN_{s_*}|^2\dd
 \leq C\frac{|B_R|}{(R-r)^2}\bigl(\sup_{B_R}Q\bigr)^{s_*}.
\end{equation}
In particular, if $B_{2r}\Subset\Omega$, then
\begin{equation}\label{eq:moment-est}
 \int_{B_r}|DN_{s_*}|^2\dd\leq C r^{-2}\int_{B_{2r}}Q^{s_*}\dd.
\end{equation}
Here and below the Hessian on $\{Q>0\}$ is the classical Hessian of the nondegenerate equation.
\end{theorem}

Observe that the conclusion extends to every $s\geq s_*$: indeed,
\begin{equation*}
 N_s=|N_{s_*}|^{s/s_*-1}N_{s_*},
\end{equation*}
and the power map on the right-hand side is locally Lipschitz (the identity map if $s=s_*$). Since $N_{s_*}$ is locally bounded, the Sobolev chain rule gives $N_s\in W^{1,2}_{\rm loc}$. In particular, the earlier sufficient range $s>m$ is contained in the range obtained here, which now also includes the interval $[s_*,m]$.

The proof does not rely only on a pointwise algebraic inequality for the Hessian. The decomposition into the normal component, the tangential derivative of the gradient magnitude, and the tangential Hessian is familiar from geometric studies of the $p$-Laplace equation; see, for example, \cite[Proposition 4.4 and (4.8)]{FMP} in the setting of capacitary potentials. In the present problem the longitudinal coefficient in the differentiated equation for $Q$ is negative when $s<m$. To estimate this term, we first use the equation on a sufficiently large interior ball to obtain a lower bound for the weighted energy of the tangential derivative of $Q$ when its magnitude has a prescribed contrast. A covering by balls chosen according to such contrasts then gives a bound for the longitudinal term in terms of quantities with positive coefficients. This permits the absorption for $s$ in a small interval below $m$. The exponent $s_*$ can be specified by a finite expression in $n,p$ alone. This quantitative choice, including the constants in the geometric and covering estimates, is given in Section~\ref{sec:effective-exponent}. The exponent is not claimed to be optimal.

The argument is robust with respect to small right-hand sides after rescaling. Let $F=\AP(Du)$, and for $\lambda>0$ define
\begin{equation}\label{eq:floor}
 Z_{s,\lambda}=F\max\{|F|,\lambda^P\}^{s/(2P)-1}.
\end{equation}

\begin{theorem}[Bounded right-hand side]\label{thm:source}
Let $n,p$ be as in Theorem \ref{thm:main}. There exist an exponent $s_*=s_*(n,p)\in(0,m)$, the same as in Theorem \ref{thm:main}, and constants $K,C$ depending only on $n,p$ such that the following holds. Let $u$ solve \eqref{eq:Equation} and let $B_r\subset B_R\Subset\Omega$ be concentric balls. Put
\begin{equation*}
 M=\sup_{B_R}Q,\qquad F_0=\norm{f}_{L^\infty(B_R)},\qquad
 \lambda=K(RF_0)^{1/P}.
\end{equation*}
If $F_0>0$, then
\begin{align}
 \int_{B_r\cap\{Q\geq\lambda\}}Q^{s_*-2}|D^2u|^2\dd
 &\leq C\frac{|B_R|}{(R-r)^2}(M^{s_*}+\lambda^{s_*}),\label{eq:source-high}\\
 \int_{B_r}|DZ_{s_*,\lambda}|^2\dd
 &\leq C\frac{|B_R|}{(R-r)^2}(M^{s_*}+\lambda^{s_*}).\label{eq:source-floor}
\end{align}
The Hessian in \eqref{eq:source-high} is the weak Hessian locally on the open set $\{Q>0\}$. If $F_0=0$, the estimates are interpreted as \eqref{eq:main-est}.
\end{theorem}

We recall that the physical flux $F$ is locally in $W^{1,2}$ for square-integrable right-hand sides by the result of Cianchi and Maz'ya \cite[Theorem 2.1]{CM}. In the power-function case of the result of Karppinen and Sarsa \cite[Theorem 1.3 and Section 3.3]{KS}, the sufficient restriction for square integrability of the differentiated untruncated field with bounded right-hand sides is $\beta\geq p-2$, which in our notation reads $s\geq2P$. Under the stronger data assumption $f\in W^{1,1}\cap L^q$, $q>n$, Montoro, Muglia and Sciunzi \cite[Theorem 1.1]{MMS} obtained the untruncated power range $s>P$ up to a smooth boundary. Our assertion concerns a different estimate, with no derivatives of $f$ and with a restriction to gradients above the source scale. For a fixed $\lambda>0$, membership of $Z_{s,\lambda}$ in $W^{1,2}_{\rm loc}$ on bounded-gradient subsets already follows from the result of Cianchi and Maz'ya and the Lipschitz chain rule. For $0<s<2P$, however, the global Lipschitz bound of the map from $F$ to $Z_{s,\lambda}$ is of order $\lambda^{s/2-P}$. Its direct use in the physical-flux estimate therefore gives, in the squared estimate, an amplitude term of order $\lambda^{s-2P}M^{2P}$, rather than $M^s$. The estimates in Theorem \ref{thm:source} control the high-gradient Hessian and the differentiated regularized flux at the source scale, without this loss when $M/\lambda$ is large. The example in Subsection \ref{sec:example} shows why the regularization in \eqref{eq:floor} is needed for bounded nonzero $f$.

The classical local $W^{2,2}$ regularity for $p<3+2/(n-2)$ is due to Manfredi and Weitsman \cite[Lemma 5.1]{MW}; see also \cite{HS} for global second-order estimates. Higher integrability of the Hessian in that strict exponent range was obtained by the Cordes argument and was reproved in \cite[Corollary 1.2]{DFZZ}. The new range in Corollary \ref{cor:unweighted-higher} includes the exponent $p=3+2/(n-2)$ and a right neighborhood of it in every dimension $n\geq3$. An explicit size of this right neighborhood is given in Corollary~\ref{eff:delta}. In dimension three, the unrefereed online manuscript of Jin \cite[Proposition~4.1]{Jin} states a weighted Hessian estimate based on the curvature of level surfaces. Assuming that estimate, restated as \eqref{eq:jin-energy} below, the higher-integrability argument of Theorem \ref{thm:anchored-native} gives the conditional extension to every finite $p>2$ in Subsection \ref{sec:three}. The present paper does not verify the proof of this external estimate. The all-dimensional results, including the unconditional three-dimensional range $2<p\leq5+\delta_3$ of Corollary~\ref{eff:delta}, do not use it. We do not address the conjectured optimal H\"older exponent of $Du$.

\section{Preliminaries and geometric estimates}

\subsection{Preliminary estimates and notation}

We use the standard local gradient regularity of scalar solutions of \eqref{eq:Equation}, including the interior estimates of DiBenedetto and Tolksdorf \cite{DiB} and \cite[Theorem 1]{Tolksdorf}. In a form sufficient here, if a weak solution on $B_l$ satisfies $\sup_{B_l}Q\leq A$, then, with $F_0=\|f\|_{L^\infty(B_l)}$, for some $\alpha_0=\alpha_0(n,p)>0$,
\begin{equation}\label{eq:holder}
 |Du(x)-Du(y)|\leq C\bigl(A+(lF_0)^{1/P}\bigr)
 \left(\frac{|x-y|}{l}\right)^{\alpha_0}
 \quad(x,y\in B_{l/2}).
\end{equation}
Indeed, put $A_0=A+(lF_0)^{1/P}$ and assume that $A_0>0$. For a ball centred at $x_0$, set
\begin{equation*}
 v(y)=\frac{u(x_0+ly)-u(x_0)}{lA_0},\qquad y\in B_1.
\end{equation*}
Then $|Dv|\leq1$, $|v|\leq1$ and
$-\diver(\AP(Dv))=lA_0^{-P}f(x_0+l\,\cdot)$, whose right-hand side has absolute value at most one. Applying the interior gradient estimate to $v$ and returning to the original variables gives \eqref{eq:holder}. If $A_0=0$, the estimate is immediate.
In the homogeneous equation we shall also use
\begin{equation}\label{eq:grad-sup}
 \sup_{B_{l/2}}Q\leq C\left(\avg_{B_l}Q^p\dd\right)^{1/p}.
\end{equation}
These estimates are not sharp H\"older estimates. In particular, no endpoint regularity is used in the proofs.

We shall use the local $W^{1,2}$ regularity of the physical flux $F$ for a square-integrable right-hand side, proved by Cianchi and Maz'ya \cite[Theorem 2.1]{CM}. Let $B\Subset\Omega$ be an interior ball and choose $f_k\in C_c^\infty(B)$ converging to $f$ in $L^2(B)$, uniformly bounded in $L^\infty(B)$. Let $u_k-u\in W_0^{1,p}(B)$ be the weak Dirichlet solution with datum $f_k$, which exists by strict convexity and coercivity of the usual energy. Taking $u_k-u$ as a test function in the difference of the two equations, by the monotonicity inequality for $\AP$ we obtain
\begin{align*}
 c_p\int_B|Du_k-Du|^p\,dx
 &\leq\int_B(\AP(Du_k)-\AP(Du))\cdot(Du_k-Du)\,dx\\
 &=\int_B(f_k-f)(u_k-u)\,dx\\
 &\leq C_B\|f_k-f\|_{L^2(B)}\|Du_k-Du\|_{L^p(B)}.
\end{align*}
Here the last inequality uses the Poincar\'e inequality and $p>2$. Dividing by the last factor when it is nonzero gives $Du_k\to Du$ in $L^p(B)$. The inequality
$|\AP(\xi)-\AP(\eta)|\leq C_p(|\xi|+|\eta|)^{p-2}|\xi-\eta|$
and H\"older's inequality imply strong convergence of $\AP(Du_k)$ to $\AP(Du)$ in $L^{p/(p-1)}(B)$. After a subsequence, $u_k$ and $Du_k$ also converge almost everywhere to $u$ and $Du$, by Poincar\'e and the strong gradient convergence. Thus $u$ belongs locally to the solution class of \cite[Theorem 2.1 and Section 5]{CM}.
On a ball where $Q$ is bounded above and bounded away from zero, the inverse of $\AP$ is smooth with bounded derivative on the relevant range. Consequently $Du\in W^{1,2}$ locally there. More quantitatively, on a ball of radius $l$ with $cA\leq Q\leq CA$ and $lF_0/A^P\leq C_0$, rescaling length by $l$ and gradient by $A$, applying the quoted local flux estimate on the normalized ball, and using the inverse map yields
\begin{equation}\label{eq:warmH}
 \int_{B_{c_1 l}}|D^2u|^2\dd\leq C A^2l^{n-2},
\end{equation}
on a fixed smaller ball, with constants depending on the prescribed comparisons and $n,p$.
To check the constant in \eqref{eq:warmH}, write
$v(y)=(u(x_0+ly)-u(x_0))/(lA)$. The normalized flux
$\widetilde F=\AP(Dv)$ is bounded, and the normalized source has bounded $L^2$ norm by the assumption on $lF_0/A^P$. The local estimate of \cite[Theorem 2.1]{CM}, applied on a fixed interior subball of $B_1$, gives $\int|D\widetilde F|^2\leq C$. On the interval of comparable positive gradients, the radial inverse of $\AP$ can be extended to a Lipschitz map by a cutoff near its origin. The weak chain rule therefore gives $\int|D^2v|^2\leq C$ there. Since $D^2u(x_0+ly)=A l^{-1}D^2v(y)$, changing variables gives \eqref{eq:warmH}.
In the homogeneous equation, on each relatively compact subset of $\wm=\{Q>0\}$ the normalized equation $\mathbb E:D^2u=0$ is uniformly elliptic and its coefficients are analytic functions of $Du$. Classical interior regularity first gives smoothness, and Morrey's analytic regularity theorem \cite[Section 5.8]{Morrey} gives real analyticity there. Uniform higher derivative estimates on strictly positive-gradient buffered balls are also available there. For the inhomogeneous equation only \eqref{eq:warmH}, not a pointwise Hessian estimate, will be used.

On $\wm$, put
\begin{equation*}
 e=Du/Q,\qquad \ell=\log Q,\qquad
 h=e\cdot D\ell,\qquad g=D\ell-he,\qquad
 \mathbb E=I+b e\otimes e.
\end{equation*}
We write $H=D^2u$ in the orthogonal decomposition $\R e\oplus e^\perp$:
\begin{equation*}
 \frac{H}{Q}=\begin{pmatrix}h&g^T\\g&S\end{pmatrix},
 \qquad S_0=S-\frac{\tr S}{n-1}I_{e^\perp}.
\end{equation*}
The equation gives
\begin{equation}\label{eq:trace}
 \tr S=-Ph-j,\qquad j=f/Q^P.
\end{equation}
These notations are used only on $\wm$; quantities with the factor $\mathbf 1_{\wm}$ are zero elsewhere. In the homogeneous problem $j=0$. The identities $D_e e=g$ and $D_\tau e=S\tau$ for $\tau\in e^\perp$ follow by differentiating $e=Du/Q$.

\subsection{Local geometric nondegeneracy}\label{sec:charge}

The key ingredient is a lower bound for the tangential derivative $g$ when the gradient magnitude has a fixed contrast on a ball with a sufficiently large filled buffer.

\begin{lemma}[Filled-ball nondegeneracy]\label{lem:charge}
Fix a compact interval $I=[s_-,s_+]\Subset(0,\infty)$. There exist $L>4$, $\delta_0>0$ and $c_0>0$, depending only on $n,p,I$, such that the following holds. Suppose $u$ solves \eqref{eq:Equation} on $B_L$, with
\begin{equation}\label{eq:contrast}
 Q(0)\leq\frac12,\qquad \max_{\overline B_1}Q=1,
 \qquad \sup_{B_L}Q\leq L^{2/s_-},\qquad \norm{f}_{L^\infty(B_L)}\leq\delta_0.
\end{equation}
Then for every $s\in I$,
\begin{equation}\label{eq:charge}
 \int_{B_2\cap\{Q\geq3/4\}}Q^s|g|^2\dd\geq c_0.
\end{equation}
\end{lemma}

\begin{proof}
We first give the exclusion in the homogeneous limiting equation. Take $L$ so large that
\begin{equation*}
 \exp\left(\frac{2(n-1)}{P(L-1)}\right)<2.
\end{equation*}
Suppose that a homogeneous solution satisfies the contrast in \eqref{eq:contrast}, but $g=0$ on a nonempty positive-gradient open patch containing a point of maximum gradient in $\overline B_1$. The field $g$ is analytic on each connected component of $\wm\cap B_L$. Hence $g=0$ on the connected positive-gradient component $W$ containing this point.

For $x\in W$, integral curves of $e$ are straight lines while they remain in $W$, since $D_e e=g=0$. The matrix $De$ vanishes on the normal direction and is symmetric on the tangent space. Along a straight normal line $x+te(x)$ the identity
\begin{equation*}
 D_e(De)=-(De)^2
\end{equation*}
holds. More precisely, with $(De)_{ij}=\partial_j e_i$, differentiating $e_k\partial_k e_i=0$ gives
\begin{equation*}
 e_k\partial_k\partial_j e_i
 =-\sum_k(\partial_k e_i)(\partial_j e_k)
 =-\bigl((De)^2\bigr)_{ij}.
\end{equation*}
Since $e$ is constant along the normal line, the tangent plane there is fixed. The symmetric restriction of $De$ to that plane solves the matrix equation $A'=-A^2$ with its initial value, whose solution is $A(t)=A(0)(I+tA(0))^{-1}$ as long as the inverse exists. Its tangential eigenvalues are therefore $\kappa_i(t)=\kappa_i(0)/(1+t\kappa_i(0))$, $1\leq i\leq n-1$, on the interval of existence before a zero denominator. The trace relation \eqref{eq:trace} with $j=0$ gives $D_e\log Q=-(\tr S)/P$. Integrating along the line yields
\begin{equation}\label{eq:normalQ}
 Q(x+te(x))=Q(x)\prod_{i=1}^{n-1}(1+t\kappa_i(0))^{-1/P}.
\end{equation}
All factors are positive on the maximal interval containing zero.

A finite exit from $W$ before a zero denominator, while the line remains in $B_L$, is impossible: the right side of \eqref{eq:normalQ} has a finite strictly positive limit, so the limiting point has a positive gradient and belongs to the same component. A zero denominator before reaching the boundary of $B_L$ is also impossible, since it would make $Q$ unbounded at an interior point by \eqref{eq:normalQ}. It follows that for every $x\in W\cap\overline B_1$, the line exists in both directions for $|t|<L-1$, and
\begin{equation*}
 |\kappa_i(0)|\leq (L-1)^{-1},\qquad
 |D\log Q(x)|=|h(x)|\leq\frac{n-1}{P(L-1)}.
\end{equation*}
Connect the maximum point to the origin by a line segment in $\overline B_1$. The bound for $D\log Q$ prevents the segment from leaving $W$ at a point where $Q=0$, and hence it applies along the entire segment. The resulting gradient ratio is strictly less than $2$, contradicting \eqref{eq:contrast}.

To prove the lemma, suppose no uniform pair $(\delta_0,c_0)$ exists for the chosen $L$. We may then take a sequence of solutions satisfying \eqref{eq:contrast}, with $\norm{f_k}_\infty\to0$, exponents $s_k\in I$, and left sides of \eqref{eq:charge} tending to zero. After subtraction of a constant and passage to a subsequence, \eqref{eq:holder} and the gradient bound yield convergence in $C^1_{\rm loc}(B_L)$ to a homogeneous solution. The contrast is preserved on $\overline B_1$. On a small neighborhood of a maximum point of the limiting gradient, the gradients of the sequence are uniformly bounded away from zero and exceed $3/4$. By \eqref{eq:warmH}, the Hessians converge weakly in $L^2$ there after passing to a subsequence. The coefficients expressing $g$ in terms of the Hessian and the normalized gradient converge uniformly. The vanishing charge implies $g=0$ on that limiting positive-gradient neighborhood. The preceding homogeneous exclusion gives the contradiction.
\end{proof}

\subsection{A weighted covering inequality}\label{sec:covering}

Choose a compact interval $I=[s_-,s_+]$ with
\begin{equation*}
 0<s_-<m<s_+<2P.
\end{equation*}
Constants in this section are uniform for $s\in I$. Let $\chi\in C^\infty(\R)$ be nondecreasing, zero on $(-\infty,1]$, one on $[2,\infty)$, and between zero and one everywhere. For $\varepsilon>0$ write $\chi_\varepsilon=\chi(Q/\varepsilon)$.

\begin{lemma}[Weighted covering estimate]\label{lem:korn}
Let $B_r\subset B_t\subset B_R\Subset\Omega$ be concentric. Suppose
\begin{equation*}
 RF_0/\varepsilon^P\leq\delta_0
\end{equation*}
for a sufficiently small constant depending only on $n,p,I$. Then for every $s\in I$,
\begin{equation}\label{eq:korn}
 \int_{B_r}\chi_\varepsilon Q^{s-2}|H|^2\dd
 \leq C_K\int_{B_t}\chi_\varepsilon Q^s|g|^2\dd
       + C_K\frac{|B_R|M^s}{(t-r)^2},
\end{equation}
where $M=\sup_{B_R}Q$. For $f=0$ there is no restriction on $\varepsilon$.
\end{lemma}

\begin{proof}
Set $d=t-r$ and $r_{\rm cap}=d/4$. For a point $x\in B_r$ with $Q(x)>\varepsilon$, let $A_x(l)=\max_{\overline B_l(x)}Q$. Define $l_x$ to be the first radius in $(0,r_{\rm cap}]$ at which $A_x(l_x)=2Q(x)$, if such a radius exists; otherwise set $l_x=r_{\rm cap}$. Thus $A_x(l_x)\leq 2Q(x)$. Smallness of $RF_0/Q(x)^P$, together with \eqref{eq:holder}, shows that there is a fixed small $\kappa>0$ such that
\begin{equation*}
 \tfrac9{10}Q(x)\leq Q\leq\tfrac{11}{10}Q(x)
 \quad\text{on } B_{20\kappa l_x}(x).
\end{equation*}
On a fixed smaller ball, \eqref{eq:warmH} hence yields
\begin{equation}\label{eq:smallball}
 \int_{B_{5\kappa l_x}(x)}\chi_\varepsilon Q^{s-2}|H|^2\dd
 \leq C\chi\bigl(11Q(x)/(10\varepsilon)\bigr)Q(x)^s l_x^{n-2}.
\end{equation}

If there is a doubling radius, choose $\rho_x\in[l_x,r_{\rm cap}]$ maximizing $A_x(l)^s/l^2$ on this interval and write $A_x=A_x(\rho_x)$. In particular,
\begin{equation}\label{eq:record}
 \frac{Q(x)^s}{l_x^2}\leq\frac{A_x^s}{\rho_x^2},\qquad
 A_x\geq2Q(x).
\end{equation}
(Using the doubled value only strengthens the first inequality.) If there is no doubling, or if $\rho_x\geq r_{\rm cap}/L$, where $L$ will be chosen large and fixed below, the bound on the right of \eqref{eq:smallball} is at most $C M^s l_x^n/d^2$. These balls are dealt with by volume packing.

Consider the remaining case $\rho_x<r_{\rm cap}/L$. The maximizing property implies
\begin{equation*}
 A_x(l\rho_x)\leq l^{2/s}A_x,\qquad 1\leq l\leq L.
\end{equation*}
The rescaled solution
\begin{equation*}
 u_x(z)=\frac{u(x+\rho_x z)-u(x)}{\rho_x A_x}
\end{equation*}
is defined on the filled ball $B_L$ and has the gradient contrast of Lemma \ref{lem:charge}. Its right-hand side has norm at most $\rho_x F_0/A_x^P\leq\delta_0$. Increasing $L$ if necessary, the preceding lemma yields
\begin{equation}\label{eq:charge-original}
 C_x=B_{2\rho_x}(x)\cap\{Q\geq3A_x/4\},\qquad
 \int_{C_x}Q^s|g|^2\dd\geq c_0 A_x^s\rho_x^{n-2}.
\end{equation}
All points of $C_x$ have gradient magnitude between $3A_x/4$ and $C A_x$, by the record property. Moreover $Q\geq3A_x/4\geq3Q(x)/2$ there, so monotonicity of $\chi$ permits the same cutoff in \eqref{eq:charge-original} as in \eqref{eq:smallball}. Combining \eqref{eq:record}--\eqref{eq:charge-original},
\begin{equation}\label{eq:pay}
 \int_{B_{5\kappa l_x}(x)}\chi_\varepsilon Q^{s-2}|H|^2\dd
 \leq C(l_x/\rho_x)^n\int_{C_x}\chi_\varepsilon Q^s|g|^2\dd.
\end{equation}

Select a countable disjoint family of the balls $B_{\kappa l_x}(x)$ whose fivefold dilations cover the centers under consideration. The additive-class cost is bounded by $C M^s|B_R|/d^2$ because these small balls are disjoint and lie in $B_t$. We next verify the weighted overlap for the remaining charge sets; their larger radii need not be disjoint.

For any $z$ in such a charge set $C_x$ define its own first upward doubling radius $l_z$ with the same spatial cap $r_{\rm cap}$. This is meaningful because $z\in B_{r+2r_{\rm cap}/L}$ and $B_{r_{\rm cap}}(z)\subset B_t$. There are constants $c,C>0$, independent of $x,z,\varepsilon$, such that
\begin{equation}\label{eq:radius-comp}
 c\rho_x\leq l_z\leq C\rho_x.
\end{equation}
For the lower bound, apply \eqref{eq:holder} to a ball of radius comparable to $\rho_x$ about $z$, using the bound on $B_{3\rho_x}(x)$ and $Q(z)\geq3A_x/4$. The gradient cannot double on a sufficiently small fixed fraction of this ball. For the upper bound choose a fixed sufficiently large number $C_1$. If $l_z>C_1\rho_x$, the bound $Q\leq2Q(z)\leq C A_x$ on $B_{C_1\rho_x}(z)$ and \eqref{eq:holder} imply that the oscillation of $Du$ between $x$ and a maximum point in $\overline B_{\rho_x}(x)$ is less than $A_x/2$ for $C_1$ large enough. This contradicts their gradient-magnitude contrast. Here the normalized source term is bounded by $((RF_0)/A_x^P)^{1/P}\leq\delta_0^{1/P}$, and the oscillation factor is $O(C_1^{-\alpha_0})$. Choose $L>C_1$ from the outset, with $L$ also large enough for the homogeneous exclusion in Lemma \ref{lem:charge}; then $C_1\rho_x<r_{\rm cap}$. The constants determining $C_1$ use only the record bound at radius $3\rho_x$, not $L^{2/s_-}$, so this choice is not circular.

At a fixed point $z$ belonging to one or more selected $C_{x_i}$, formula \eqref{eq:radius-comp} makes all the associated $\rho_{x_i}$ comparable to the same $l_z$. Their centers lie in a common ball of radius $C l_z$. The disjoint volumes of $B_{\kappa l_{x_i}}(x_i)$ therefore give
\begin{equation*}
 \sum_{i:z\in C_{x_i}}(l_{x_i}/\rho_{x_i})^n\leq C.
\end{equation*}
Summing \eqref{eq:pay} and integrating this bound gives the charge-side term of \eqref{eq:korn}. The covering includes the support of $\chi_\varepsilon$ in $B_r$. This proves the lemma.
\end{proof}

For the proofs in the next two sections, fix $0<\delta_0\leq1$ below all the smallness thresholds used in Lemmas~\ref{lem:charge} and \ref{lem:korn}, including the fixed-scale gradient comparisons in the covering proof. Decreasing this common threshold preserves all their conclusions.

\section{Proof of the weighted second-order estimates}

\subsection{The differentiated magnitude equation and absorption}\label{sec:absorption}

We give the distributional calculation with a bounded right-hand side, so that the homogeneous and inhomogeneous proofs use the same covering and coefficient choice.

\begin{lemma}[Magnitude identity on the positive-gradient set]\label{lem:bochner}
On $\wm$, for $s>0$, one has
\begin{align}
 \diver(Q^s\mathbb E D\ell)
 &=Q^s\bigl(|S|^2+s|g|^2+P(s-P)h^2\bigr)
     -Q^{s-P}e\cdot Df\label{eq:raw-bochner}\\
 &=Q^s\left(|S_0|^2+s|g|^2+P(s-m)h^2
       +\frac{2P}{n-1}jh+\frac{j^2}{n-1}\right)
       -Q^{s-P}e\cdot Df. \label{eq:bochner}
\end{align}
Here the term involving $Df$ is interpreted by distributional integration by parts against compactly supported tests supported in the positive-gradient set.
\end{lemma}

\begin{proof}
On a compact subset of $\wm$, put $B=Q^{p-2}\mathbb E$. Since $Du\in W^{1,2}$ locally on this set, differentiating the weak equation in the $i$-th coordinate gives
\begin{equation*}
 \int BDu_i\cdot D\zeta\,dx=-\int f\,\partial_i\zeta\,dx .
\end{equation*}
This formula holds also for compactly supported $W^{1,2}$ tests there, by density, because $B$ is locally bounded, $Du_i\in L^2$, and $f\in L^2$. Let $\varphi$ be a compactly supported smooth function on this set. Taking $\zeta=u_i\varphi$ and summing over $i$, we obtain
\begin{equation*}
 \int\varphi\tr(HBH)\,dx
 +\int BHDu\cdot D\varphi\,dx
 =-\int f(\varphi\Delta u+Du\cdot D\varphi)\,dx.
\end{equation*}
Observe that $HDu=QDQ=Q^2D\ell$, and hence $BHDu=Q^p\mathbb E D\ell$. It follows that
\begin{equation*}
 \diver(Q^p\mathbb E D\ell)
 = Q^{p-2}\tr(H\mathbb E H)-Du\cdot Df
\end{equation*}
on this set as a distribution. Here the last term is defined by $\langle -Du\cdot Df,\varphi\rangle=\int f\diver(Du\varphi)\,dx$.

Multiplication on this set by $Q^{s-p}$, using the usual weak product rule, subtracts $(p-s)Q^s\mathbb E D\ell\cdot D\ell$ from the quadratic term. In the frame $e\oplus e^\perp$,
\begin{equation*}
 Q^{-2}\tr(H\mathbb E H)=P h^2+(P+1)|g|^2+|S|^2,\qquad
 \mathbb E D\ell\cdot D\ell=P h^2+|g|^2.
\end{equation*}
This proves \eqref{eq:raw-bochner}. Finally \eqref{eq:trace} implies
$|S|^2=|S_0|^2+(Ph+j)^2/(n-1)$; substitute this relation and use the definition of $m$. All weak products may be obtained by mollification on compact positive-gradient subsets, where powers and inverses of $Q$ have bounded derivatives. The continuous Sobolev coefficients and test functions are approximated both uniformly and in $W^{1,2}$; this justifies their pairing with the $L^1$ quadratic Hessian densities, as well as the terms paired with $Df$.
\end{proof}

\begin{lemma}[Truncated production estimate]\label{lem:production}
Let $I$ and $\chi_\varepsilon$ be as in Lemma \ref{lem:korn}. For $s\in I$, $s\leq m$, and any $\eta\in C_c^\infty(B_R)$, $0\leq\eta\leq1$, one has, for each $\zeta>0$,
\begin{align}
 \int\eta^2\chi_\varepsilon Q^s(|S_0|^2+s|g|^2)\dd
 &\leq C\bigl(P(m-s)+\zeta\bigr)
              \int\eta^2\chi_\varepsilon Q^s h^2\dd \notag\\
 &\quad+C_\zeta\int Q^s|D\eta|^2\dd
          +C_\zeta F_0^2\varepsilon^{s-2P}|\supp\eta|.\label{eq:production}
\end{align}
The constants can be chosen uniform for $s\in I$. No last term is needed when $F_0=0$.
\end{lemma}

\begin{proof}
Test \eqref{eq:bochner} with $\eta^2\chi_\varepsilon$. This test is zero in a neighborhood of the critical set on its compact support and is admissible by Sobolev approximation on the positive-gradient set. The term arising from the derivative of $\chi_\varepsilon$ on the divergence side is, with favorable sign,
\begin{equation*}
 \mathcal G=\int\eta^2\chi'(Q/\varepsilon)\frac{Q^{s+1}}{\varepsilon}
                       (Ph^2+|g|^2)\dd\geq0.
\end{equation*}
Integrate the $Df$ term by parts. In particular,
\begin{equation*}
 \diver(Q^{s-P}e)=Q^{s-P}\bigl((s-2P)h-j\bigr)
\end{equation*}
by \eqref{eq:trace}. By integration by parts of the source term and the trace relation, we obtain
\begin{align*}
 &\int\eta^2\chi_\varepsilon Q^s(|S_0|^2+s|g|^2)\,dx+\mathcal G\\
 &=P(m-s)\int\eta^2\chi_\varepsilon Q^s h^2\,dx
 -2\int\eta\chi_\varepsilon Q^s\mathbb E D\ell\cdot D\eta\,dx\\
 &\quad -(s-2m)\int\eta^2\chi_\varepsilon Q^s jh\,dx
 +\frac mP\int\eta^2\chi_\varepsilon Q^s j^2\,dx\\
 &\quad -2\int\eta\chi_\varepsilon Q^s j e\cdot D\eta\,dx
 -\int\eta^2\chi'(Q/\varepsilon)\frac{Q^{s+1}}\varepsilon jh\,dx.
\end{align*}
For the second term on the right, write $\mathbb E D\ell=Ph e+g$. By Young's inequality, its absolute value is at most
\begin{equation*}
 \zeta\int\eta^2\chi_\varepsilon Q^s h^2\,dx
 +\frac{s_-}{4}\int\eta^2\chi_\varepsilon Q^s|g|^2\,dx
 +C_\zeta\int Q^s|D\eta|^2\,dx .
\end{equation*}
The next $jh$ term is bounded by another $\zeta$ times the longitudinal integral and $C_\zeta\int\eta^2\chi_\varepsilon Q^s j^2\,dx$. For the final term, one has
\begin{equation*}
 |jh|\leq\tfrac P2h^2+\tfrac1{2P}j^2.
\end{equation*}
Multiplying this inequality by the nonnegative cutoff weight, the first contribution is at most $\mathcal G/2$. Since $\chi'$ is supported where $1<Q/\varepsilon<2$, the second is at most $C F_0^2\varepsilon^{s-2P}|\supp\eta|$. All remaining source terms are bounded by a constant times this expression and the ordinary cutoff integral. Combining the inequalities and using $s\geq s_->0$ yields \eqref{eq:production}.
\end{proof}

\begin{proof}[Proof of Theorems \ref{thm:main} and \ref{thm:source}, high-gradient estimate]
Fix $I$ about $m$ as in Lemma \ref{lem:korn}, and use its uniform constant $C_K$. For $s\in I$ with $s\leq m$ set
\begin{equation*}
 X_\varepsilon(l)=\int_{B_l}\chi_\varepsilon Q^{s-2}|H|^2\dd.
\end{equation*}
Apply \eqref{eq:korn} to $B_r\subset B_{(r+t)/2}$, and \eqref{eq:production} with a cutoff equal to one on $B_{(r+t)/2}$, supported in $B_t$ and with $|D\eta|\leq C/(t-r)$. Because $Q^sh^2\leq Q^{s-2}|H|^2$, these two inequalities imply
\begin{equation*}
 X_\varepsilon(r)\leq C_1\bigl(P(m-s)+\zeta\bigr)X_\varepsilon(t)
 + C_\zeta |B_R|\left(\frac{M^s}{(t-r)^2}
                  +F_0^2\varepsilon^{s-2P}\right).
\end{equation*}
All constants determining the first coefficient are uniform for $s$ in the fixed interval. Choose $\zeta>0$ small and then choose $s_*=m-\sigma\in I$ with $\sigma>0$ small enough that the first coefficient is less than $1/4$. This choice depends only on $n,p$.

Under $RF_0/\varepsilon^P\leq\delta_0$,
$F_0^2\varepsilon^{s_*-2P}\leq C\varepsilon^{s_*}/R^2$. Thus ordinary geometric iteration in radii increasing from $r$ to $R$ gives
\begin{equation}\label{eq:truncated-est}
 \int_{B_r}\chi_\varepsilon Q^{s_*-2}|H|^2\dd
 \leq C\frac{|B_R|}{(R-r)^2}(M^{s_*}+\varepsilon^{s_*}),
\end{equation}
for every admitted $\varepsilon>0$. To see the iteration directly, take $r_j=r+(R-r)(1-2^{-j})$ and sum the resulting geometric series with ratio less than one after the factor $4^j$. The terminal integral stays finite for fixed $\varepsilon$, by local $W^{1,2}$ regularity of $Du$ on the compact subset $\{Q\geq\varepsilon\}$ of the positive-gradient set; no untruncated integral is assumed.

In the homogeneous equation send $\varepsilon\downarrow0$ by monotone convergence. This gives the weighted Hessian estimate of Theorem \ref{thm:main}. The vector fields
\begin{equation*}
 N_{s_*,\varepsilon}=(Q^{s_*/2}-(2\varepsilon)^{s_*/2})_+e
\end{equation*}
have weak derivatives controlled in square by $C Q^{s_*-2}|H|^2\mathbf 1_{\{Q>2\varepsilon\}}$. To justify the global derivative on a compact ball, they vanish on the closed set $Q\leq2\varepsilon$ and are given by a Lipschitz chain rule on a buffered positive-gradient set. Equivalently they are Lipschitz compositions of the locally Sobolev flux $F$ for each fixed $\varepsilon$. They converge uniformly on compact sets to $N_{s_*}$, and their derivative norms are uniformly bounded. Weak compactness yields the asserted $W^{1,2}$ conclusion and estimate \eqref{eq:main-est}.

For \eqref{eq:moment-est}, apply \eqref{eq:main-est} on $B_r\subset B_{3r/2}$ and the bound $\sup_{B_{3r/2}}Q^{s_*}\leq C\avg_{B_{2r}}Q^{s_*}$. The latter follows from \eqref{eq:grad-sup} for nested subballs and interpolation: write $\int Q^p\leq M^{p-s_*}\int Q^{s_*}$ and apply Young's inequality to absorb an arbitrarily small multiple of the larger-ball supremum; iterate over a geometric sequence of intermediate radii. For clarity, the resulting inequality has the form $M_a\leq\theta M_b+C_\theta(b-a)^{-n/s_*}(\int_{B_{2r}}Q^{s_*})^{1/s_*}$ for $3r/2\leq a<b<2r$; choose $\theta<2^{-n/s_*}$ before iteration. This proves \eqref{eq:moment-est}.

When $F_0>0$, choose $\varepsilon=(RF_0/\delta_0)^{1/P}$ and $\lambda=2\varepsilon$. Estimate \eqref{eq:truncated-est} immediately implies \eqref{eq:source-high}. Thus $K=2\delta_0^{-1/P}$ in Theorem \ref{thm:source}, with the common threshold fixed after Lemma~\ref{lem:korn}.
\end{proof}

\subsection{Regularity at the source scale}\label{sec:sourcefloor}

We now prove the full regularized estimate \eqref{eq:source-floor} using the ordinary div--curl identity. This step is different from extending the singular weight $Q^{s_*-2}$ through a nonzero source.

\begin{proof}[Proof of Theorem \ref{thm:source}, estimate \eqref{eq:source-floor}]
Assume $F_0>0$ and use $s=s_*$ and the fixed $\lambda$ of the preceding proof. Write the vector field in \eqref{eq:floor} on the positive-gradient set as
\begin{equation*}
 Z=\rho_s(Q)e,\qquad \rho_s(Q)=Q^P\max(Q^P,\lambda^P)^{s/(2P)-1}.
\end{equation*}
It belongs to $W^{1,2}_{\rm loc}$, since it is a globally Lipschitz composition of the physical flux for each fixed positive $\lambda$ (recall $0<s<2P$). The preliminary flux regularity cited in Section 2 therefore justifies using the ordinary $L^2$ div--curl identity below.

Write $a=Q\rho_s'(Q)/\rho_s(Q)$ where differentiable; $a=P$ for $Q<\lambda$ and $a=s/2$ for $Q>\lambda$. At $Q=\lambda$, $DQ=0$ almost everywhere on the level by Sobolev locality, and either value for $a$ gives the same derivative there. At a differentiability point with $Q>0$ the matrix of $DZ$ in the adapted frame is
\begin{equation*}
 DZ=\rho_s(Q)\begin{pmatrix}ah & ag^T\\g&S\end{pmatrix}.
\end{equation*}
Set $|dZ|^2=\sum_{i<j}(\partial_i Z_j-\partial_j Z_i)^2$. The trace relation \eqref{eq:trace} and direct expansion give
\begin{align}
 |DZ|^2-|\diver Z|^2-|dZ|^2
 &=\rho_s(Q)^2\left[
   |S_0|^2+2a|g|^2+P(2a-m)h^2\right.\notag\\
 &\hspace{38mm}\left.+(2a-2m)jh-\frac{m}{P}j^2\right].\label{eq:hodge-density}
\end{align}
For example, the curl contributes $\rho_s^2(a-1)^2|g|^2$, the divergence contributes $\rho_s^2((a-P)h-j)^2$, and the squared trace part of $S$ is $(Ph+j)^2/(n-1)$. These identities prove \eqref{eq:hodge-density}.

On $0<Q<\lambda$, one has $a=P$ and $2P>m$. Observe that, by $\tr S=-Ph-j$, one has
\begin{equation*}
 |DZ|^2\leq C\rho_s^2\bigl(h^2+|g|^2+|S_0|^2+j^2\bigr)
\end{equation*}
on $0<Q<\lambda$. On this region $a=P$, so the positive coefficients $2P$ and $P(2P-m)$ in \eqref{eq:hodge-density} allow $|(2P-2m)jh|$ to be absorbed by a fixed fraction of the longitudinal term plus $Cj^2$. Thus the asserted lower bound follows, using $\rho_s^2j^2=\lambda^{s-2P}f^2$.
In particular,
\begin{equation*}
 |DZ|^2-|\diver Z|^2-|dZ|^2\geq c|DZ|^2-C\lambda^{s-2P} F_0^2.
\end{equation*}
Here we used $\rho_s(Q)^2 j^2=\lambda^{s-2P}f^2$ on that region. On the critical set itself $DZ=0$ almost everywhere by Sobolev locality. On $Q\geq\lambda$, the absolute value of the left side of \eqref{eq:hodge-density} is bounded by $C|DZ|^2$, and $|DZ|^2\leq C Q^{s-2}|H|^2$ there; this region has already been estimated.

For $\eta\in C_c^\infty(B_R)$ apply the standard identity
\begin{equation*}
 \int |D(\eta Z)|^2\dd=\int |\diver(\eta Z)|^2\dd+\int |d(\eta Z)|^2\dd.
\end{equation*}
It follows by integration by parts for smooth compact vector fields and by $W^{1,2}$ density in this class. To see the expansion, write $D(\eta Z)=\eta DZ+Z\otimes D\eta$, $\diver(\eta Z)=\eta\diver Z+Z\cdot D\eta$, and $d(\eta Z)=\eta dZ+D\eta\wedge Z$. By
\begin{equation*}
 |Z\otimes D\eta|^2=(Z\cdot D\eta)^2+|D\eta\wedge Z|^2,
\end{equation*}
the pure quadratic cutoff terms cancel. Estimating the mixed terms, which are bounded by $C\eta|DZ||Z||D\eta|$, by Young's inequality gives
\begin{equation*}
 \left|\int \eta^2 (|DZ|^2-|\diver Z|^2-|dZ|^2)\dd\right|
 \leq \zeta\int\eta^2|DZ|^2\dd + C_\zeta\int |Z|^2|D\eta|^2\dd.
\end{equation*}
Split the left side using the preceding low- and high-gradient estimates and take $\zeta$ sufficiently small. Choose $\eta=1$ on $B_r$ and supported in $B_{(r+R)/2}$, with $|D\eta|\leq C/(R-r)$. The high-gradient integral on this intermediate ball is bounded by \eqref{eq:source-high}, while $|Z|^2\leq M^s+\lambda^s$ and $F_0^2\lambda^{s-2P}\leq C\lambda^s/R^2$. This yields \eqref{eq:source-floor}.
\end{proof}

\subsection{The low-gradient obstruction for a nonzero source}\label{sec:example}

The following example, also discussed in \cite[Section 3.3]{KS}, shows that the distinction between the homogeneous theorem and \eqref{eq:source-floor} is necessary. Consider a one-dimensional solution extended to $\R^n$. Let
\begin{equation*}
 u(x)=-\frac{P}{p}|x_1|^{p/P}.
\end{equation*}
Then $Du=-\operatorname{sgn}(x_1)|x_1|^{1/P}e_1$, $F=-x_1e_1$, and $-\diver F=1$ in distributions. On $x_1\ne0$,
\begin{equation*}
 |Du|^{s-2}|D^2u|^2=\frac1{P^2}|x_1|^{s/P-2}.
\end{equation*}
The latter quantity is locally integrable across $x_1=0$ if and only if $s>P$. In particular the unregularized below-$m$ estimate of Theorem \ref{thm:main} fails for bounded nonzero $f$, even for a constant right-hand side. The regularization in Theorem \ref{thm:source} respects the natural source scale without making such an untruncated assertion.

\section{Higher integrability and distributional identities}

For $s>0$ write $N_s=Q^{s/2}e$ on $\{Q>0\}$ and $N_s=0$ on
$\{Q=0\}$. On the positive-gradient set the derivative of the radial power map on $Du$
has radial eigenvalue $(s/2)Q^{s/2-1}$ and tangential eigenvalues
$Q^{s/2-1}$. Consequently,
\begin{equation*}
 c_s Q^{s-2}|D^2u|^2\leq |DN_s|^2
 \leq C_s Q^{s-2}|D^2u|^2 .
\end{equation*}
For the exponent $s=s_*$ of the preceding local weighted Hessian
theorem, this field belongs to $W^{1,2}_{\mathrm{loc}}$ also across its
zero set, with $DN_s=0$ a.e. on that set. The following estimate improves the unanchored amplitude bound by subtracting the mean of this nonlinear field. For the earlier power range, Sarsa \cite[Theorem 1.1, Section 4 and Corollary 1.2]{Sarsa} established the corresponding arbitrary-constant estimate and higher integrability, using a division into degenerate and nondegenerate cases; see also the nonlinear-gradient estimates of Avelin, Kuusi and Mingione \cite[Section 4]{AKM}. The analogous conclusions in the Orlicz setting are given in \cite[Theorem 1.1 and Corollary 1.2]{KS}. For clarity, we formulate the transfer from the unanchored estimate separately from the power restriction. We include its proof because we will apply it to the exponent $s_*<m$ and, in dimension three, to $s=1$.

\begin{theorem}[Anchored estimate and higher integrability]\label{thm:anchored-native}
Let $n\geq3$, $p>2$, and let $u$ be a scalar $p$-harmonic function on an open set $\Omega\subset\mathbb R^n$. Suppose that for some $0<s<p$ there is a constant $C_0<\infty$ such that $N_s\in W^{1,2}_{\mathrm{loc}}(\Omega)$ and
\begin{equation}\label{eq:unanchored-input}
 \int_{B_r}|DN_s|^2\,dx\leq \frac{C_0}{r^2}\int_{B_{2r}}Q^s\,dx
 \qquad (B_{2r}\Subset\Omega).
\end{equation}
Then there are $C<\infty$ and $\varepsilon_{\rm hi}>0$, depending only on $n,p,s,C_0$, such that, for every $B_{4r}\Subset\Omega$,
\begin{equation*}
 \int_{B_r}|DN_s|^2\,dx\leq
 \frac{C}{r^2}\int_{B_{4r}}|N_s-(N_s)_{B_{4r}}|^2\,dx .
\end{equation*}
In particular $N_s\in W^{1,2+\varepsilon_{\rm hi}}_{\mathrm{loc}}(\Omega)$, and $Q^{s-2}|D^2u|^2$ on $\{Q>0\}$, extended by zero to the critical set, belongs to $L^{1+\varepsilon_{\rm hi}/2}_{\mathrm{loc}}(\Omega)$. All these conclusions apply unconditionally to the exponent $s=s_*(n,p)$ of Theorem \ref{thm:main}.
\end{theorem}

\begin{proof}
We first recall a consequence of the ordinary interior gradient estimate
for scalar $p$-harmonic functions which does not require an endpoint
H\"older estimate. For each fixed $0<s<p$,
\begin{equation*}
 \sup_{B_{3r}}Q\leq C(n,p,s)
 \left(\fint_{B_{4r}}Q^s\,dx\right)^{1/s}.
\end{equation*}
Indeed, after rescaling to a unit-size ball, the usual local gradient
bound by the $p$-energy on nested balls of radii $a<b$ gives
\begin{equation*}
 M(a)\leq C(b-a)^{-n/p}M(b)^{1-s/p}
 \left(\int_{B_4}Q^s\,dx\right)^{1/p},
 \qquad M(a)=\sup_{B_a}Q .
\end{equation*}
Young's inequality therefore gives
$M(a)\leq\theta M(b)+C_{\theta}(b-a)^{-n/s}
(\int_{B_4}Q^s)^{1/s}$. Choose a fixed $\theta<2^{-n/s}$
and iterate on radii increasing from $3$ to $4$ with dyadic gaps.
The terminal suprema are finite on the compact outer ball, so this gives
the asserted bound. In particular it applies to $s=s_*<m<p$.

Put
\begin{equation*}
 a=(N_s)_{B_{4r}},\qquad D=\fint_{B_{4r}}|N_s-a|^2\,dx .
\end{equation*}
Fix a sufficiently small $\delta=\delta(n,p,s)>0$.
If $D\geq\delta|a|^2$, then
$\fint_{B_{4r}}Q^s=D+|a|^2\leq(1+\delta^{-1})D$.
The assumed unanchored estimate \eqref{eq:unanchored-input} between $B_r$ and $B_{2r}$ gives directly
\begin{equation*}
 \int_{B_r}|DN_s|^2\,dx\leq Cr^{n-2}D .
\end{equation*}
The zero-field case is included.

Suppose instead that $D<\delta|a|^2$, with $a\neq0$, and set
$A=|a|^{2/s}$. Normalize length by $r$ and gradient amplitude by $A$,
subtracting only a constant from $u$. These operations preserve
$p$-harmonicity. The normalized fields have uniformly bounded gradient on
$B_3$, by the local supremum estimate and the identity
$\fint |N_s|^2=D+|a|^2$. The ordinary interior $C^{1,\alpha_0}$
estimate gives a common modulus of continuity for the gradients on
$B_{5/2}$, and hence, by continuity of the radial power map on bounded
sets, for the normalized $N_s$ as well. Small normalized $L^2$ variance
then gives
\begin{equation*}
 |N_s-a|\leq\tfrac12|a|\quad\text{on }B_{2r}
\end{equation*}
for sufficiently small fixed $\delta$. To see this directly, a violation
would give, by the common modulus, a ball of uniform positive normalized
radius inside $B_{5/2}$ on which $|N_s-a|/|a|\geq1/4$, contradicting
$D/|a|^2\to0$.

Thus $Q\asymp A$ on $B_{2r}$, with constants depending only on $s$.
The equation is nondegenerate throughout this ball, so its differentiated
equations there are
$\operatorname{div}(Q^{p-2}\mathbb E D u_i)=0$.
Choose the constant vector $\xi=Aa/|a|$, for which
$|\xi|^{s/2-1}\xi=a$. Let $\eta\in C_c^\infty(B_{2r})$ satisfy $0\leq\eta\leq1$, $\eta=1$ on $B_r$, and $|D\eta|\leq C/r$. Taking $\eta^2(u_i-\xi_i)$ as a test function in the $i$-th differentiated equation and summing over $i$ yields
\begin{equation*}
 \sum_i\int\eta^2 Q^{p-2}\mathbb E Du_i\cdot Du_i\,dx
 =-2\sum_i\int\eta(u_i-\xi_i)Q^{p-2}\mathbb E Du_i\cdot D\eta\,dx.
\end{equation*}
As $Q\asymp A$, uniform ellipticity and Young's inequality give
\begin{equation*}
 A^{p-2}\int\eta^2|D^2u|^2\,dx
 \leq C A^{p-2}\int|Du-\xi|^2|D\eta|^2\,dx,
\end{equation*}
Consequently,
\begin{equation*}
 \int_{B_r}|D^2u|^2\,dx
 \leq\frac{C}{r^2}\int_{B_{2r}}|Du-\xi|^2\,dx .
\end{equation*}
Here no affine function has been subtracted from the original nonlinear
equation: $\xi_i$ is a constant in a linear differentiated equation.
For $Q\asymp A$ the elementary power-map inequality gives
$|N_s-a|\asymp A^{s/2-1}|Du-\xi|$. For example, both squared
distances split into the squared radial difference and the
nonnegative angular term proportional to $1-\cos\vartheta$,
where $\vartheta$ is the angle between the vectors; the radial
powers are comparable by the mean value theorem on comparable
positive magnitudes. Multiply the last energy estimate by
$A^{s-2}$ and use the derivative comparison to obtain the claimed
anchored estimate also in the small-variance case.

Finally, put $t=2n/(n+2)\in(1,2)$. Sobolev--Poincar\'e on $B_{4r}$
and the anchored estimate imply the weak reverse H\"older inequality
\begin{equation*}
 \fint_{B_r}|DN_s|^2\,dx
 \leq C\left(\fint_{B_{4r}}|DN_s|^t\,dx\right)^{2/t}.
\end{equation*}
The local Gehring lemma with a fixed dilation (see, e.g., \cite{Giaquinta}) yields
$DN_s\in L^{2+\varepsilon_{\rm hi}}_{\mathrm{loc}}$ for some
$\varepsilon_{\rm hi}=\varepsilon_{\rm hi}(n,p,s,C_0)>0$. The weighted Hessian assertion follows from the positive-gradient derivative comparison and Sobolev locality on the critical set. Finally, for $s=s_*$ the assumptions are exactly Theorem \ref{thm:main} and \eqref{eq:moment-est}. Thus the anchored estimate and the reverse H\"older argument apply at the power $s_*<m$.
\end{proof}

\begin{corollary}[Unweighted Hessian beyond the usual exponent range]
\label{cor:unweighted-higher}
For every $n\geq3$ there is $\delta_n>0$ such that each scalar
$p$-harmonic function belongs to $W^{2,2+\varepsilon_{\rm hi}}_{\mathrm{loc}}$ for
\begin{equation*}
  2<p\leq 3+\frac{2}{n-2}+\delta_n,
\end{equation*}
where $\varepsilon_{\rm hi}>0$ can depend on $n,p$.
\end{corollary}

\begin{proof}
Write $p_c=3+2/(n-2)$; equivalently $m(p_c)=2$.
We explain why the choice of a positive gain below $m(p)$ in the
homogeneous estimate can be made uniformly for $p$ in a sufficiently
small closed interval $J$ about $p_c$. First fix $J\Subset(2,\infty)$
small enough that a single compact interval $I=[s_-,s_+]$ satisfies
$0<s_-<m(p)<s_+<2(p-1)$ for all $p\in J$, with the two strict gaps
uniformly positive.

Write $J=[p_-,p_+]$. Lemma~\ref{eff:holder-lemma}, proved independently in Section~\ref{sec:effective-exponent}, gives a common positive H\"older exponent by evaluating \eqref{eff:alpha} at $p_+$, with the same numerical prefactor in \eqref{eff:pair-holder}. After length and gradient normalization, this is the uniform homogeneous version of \eqref{eq:holder} needed here.

The strictly positive-gradient Hessian bounds are uniform as well. On a normalized buffered ball where $0<q_0\leq Q\leq A$, with $q_0<1\leq A$, the differentiated coefficient matrices have the common ellipticity bounds
\begin{equation*}
 \lambda_{\rm ell}=q_0^{p_+-2},\qquad
 \Lambda_{\rm ell}=(p_+-1)A^{p_+-2}.
\end{equation*}
Use the common gradient H\"older bound just obtained and evaluate $M_c$ in \eqref{eff:warm-def} at $p_+$. The coefficient estimates in the proof of Lemma~\ref{eff:warm} are then uniform on $J$. Applying Lemma~\ref{eff:linear} with these common inputs, as in the first step of that proof, gives a common Hessian bound on a smaller ball. The proofs of Lemmas~\ref{eff:holder-lemma}--\ref{eff:warm} use neither this corollary nor the higher-integrability argument of this section, so these forward references are noncircular.
In Lemma~\ref{lem:charge}, choose $L$ using
$\min_{p\in J}(p-1)>0$. Failure of a uniform homogeneous lower charge
would yield a sequence $p_k\to p_\infty\in J$ and a locally
$C^1$-convergent sequence of normalized $p_k$-harmonic functions.
The uniform gradient bound and the joint continuity of
$(p,\xi)\mapsto |\xi|^{p-2}\xi$ on bounded sets with $p\in J$ show that
the limit is $p_\infty$-harmonic; convergence of the weak equations
is thus direct. The Hessians on the uniform positive-gradient
patch at a maximum point are bounded in $L^2$, and the vanishing
charge gives $g=0$ there for the limit, exactly as in the proof of
the lemma. The normal-line argument then yields the same contradiction,
since the choice of $L$ was uniform over $J$. Consequently the charge
constant can be chosen uniformly over $p\in J$ and $s\in I$.

The constants in Lemma~\ref{lem:korn} now have the same uniformity:
all the local H\"older, positive-gradient Hessian, charge and overlap constants
just discussed are uniform. The homogeneous part of
Lemma~\ref{lem:production} uses only algebraic coefficients continuous
in $(p,s)$ and the positive lower endpoint of $I$. Its constants,
for fixed Young parameter, are therefore uniform on $J\times I$.
The absorption and radius iteration choose a single $\sigma>0$,
small enough that $s(p)=m(p)-\sigma\in I$ and the iteration coefficient
is less than $1/4$ for each $p\in J$. This proves the asserted
uniform gain, using only the homogeneous equation.

For $2<p\leq p_c$, one already has $s_*(n,p)<m(p)\leq2$.
For $p\in[p_c,p_c+\delta_n]$, choose $\delta_n>0$ so that this interval
lies in $J$ and
$m(p)\leq2+\sigma/2$. The uniform choice above gives
$s(p)\leq2-\sigma/2$. In either case the anchored estimate and
Gehring argument, applied to the corresponding exponent $s\in(0,2)$,
give $N_s\in W^{1,2+\varepsilon_{\rm hi}}_{\mathrm{loc}}$. Finally
\begin{equation*}
   Du=|N_s|^{2/s-1}N_s.
\end{equation*}
The map on the right is locally Lipschitz, including at zero, because
$2/s>1$, and $N_s$ is locally bounded. The Sobolev chain rule implies
$Du\in W^{1,2+\varepsilon_{\rm hi}}_{\mathrm{loc}}$, which proves the claim.
\end{proof}

\begin{theorem}[Distributional speed identities through the critical set]
\label{thm:cold-safe-speed}
Let $u$ be scalar $p$-harmonic as in Theorem \ref{thm:main}.
For every fixed $t\geq s_*$, extend the fields and densities below by
zero on $\{Q=0\}$. Then
\begin{equation*}
 \operatorname{div}\bigl(Q^t \mathbb E D\ell\bigr)
 = Q^t\bigl(|S_0|^2+t|g|^2+P(t-m)h^2\bigr)
 \quad\text{in }\mathcal D'(\Omega),
\end{equation*}
and
\begin{equation*}
 \operatorname{div}(Q^t e)=(t-P)Q^t h
 \quad\text{in }\mathcal D'(\Omega).
\end{equation*}
In particular the first identity at $t=m$ has a nonnegative locally
integrable density and no additional measure supported on the critical
set.
\end{theorem}

\begin{proof}
On compact positive-gradient subsets, $u$ is smooth and differentiation
of the homogeneous equation gives
\begin{equation*}
 \operatorname{div}(Q^p \mathbb E D\ell)
 =Q^{p-2}\operatorname{tr}(D^2u \mathbb E D^2u).
\end{equation*}
Multiplication by $Q^{t-p}$ and the product rule give
\begin{equation*}
 \operatorname{div}(Q^t \mathbb E D\ell)
 =Q^t\bigl(|S|^2+t|g|^2+P(t-P)h^2\bigr),
 \qquad S=D^2u|_{e^\perp}/Q.
\end{equation*}
Since $\operatorname{tr}S=-Ph$, its squared norm is
$|S_0|^2+P^2h^2/(n-1)$, which gives the first nondegenerate identity.
The second follows from $\operatorname{div}e=-Ph$ and
$D_eQ^t=tQ^t h$.

We verify the first identity across $\{Q=0\}$ without differentiating
any field there. Write $s=s_*$ and let $\chi$ be a smooth nondecreasing
function, equal to zero on $[0,1]$ and to one on $[2,\infty)$.
For $\varphi\in C_c^\infty(\Omega)$, test the nondegenerate identity by
$\chi(Q/\varepsilon)\varphi$ with $\varepsilon>0$. Its support is compactly separated
from the critical set, so this is an ordinary admissible smooth
test there. The established weighted estimate yields
\begin{equation*}
 \int_{K\cap\{Q>0\}} Q^s(h^2+|g|^2+|S|^2)\,dx<\infty,
 \qquad K=\operatorname{supp}\varphi .
\end{equation*}
Local boundedness of $Q$ and $t\geq s$ imply that every density on the
right has an integrable absolute majorant. The current on the left also
belongs to $L^1_{\mathrm{loc}}$: its norm is bounded by
$C Q^{t-s/2}(Q^{s/2}|D\ell|)$, so Cauchy--Schwarz applies locally.
The only new integral from the cutoff is bounded by
\begin{equation*}
 C\|\varphi\|_\infty
 \int_{K\cap\{\varepsilon<Q<2\varepsilon\}}
 Q^t(P h^2+|g|^2)\,dx ,
\end{equation*}
since $Q\chi'(Q/\varepsilon)/\varepsilon$ is bounded on this band.
It tends to zero as $\varepsilon\downarrow0$ by absolute continuity of
the integrable weighted energy, also when $t=s$. Dominated convergence
in the remaining terms proves the first distributional identity,
including the border exponent $t=s$.

For the second identity, alternatively use the globally Sobolev field
$Q^t e=|N_s|^{2t/s-1}N_s$. The power map here is locally Lipschitz
on bounded sets because $2t/s\geq2$. The field is locally in
$W^{1,2}$, and its weak derivative is zero a.e. on its zero set.
The nondegenerate chain computation thus gives its divergence a.e.,
hence distributionally, with no additional defect.
\end{proof}

\begin{remark}
The exact signed identity for $t<m$ is not a nonnegative
production formula. The negative longitudinal term is essential.
These estimates neither prove a pointwise second derivative at a
critical point nor give the optimal gradient H\"older conjecture.
\end{remark}

\subsection{A conditional consequence in dimension three}\label{sec:three}

This subsection is separate from the proof of the all-dimensional estimates. Its additional input is the weighted estimate stated in the unrefereed online manuscript of Jin \cite[Proposition~4.1]{Jin}. We do not verify its proof here and treat it explicitly as an external hypothesis. We state the precise hypothesis in our notation to distinguish it from the transfer argument proved below.

For a fixed finite $p\geq5$, assume that there exists $L_p<\infty$, depending only on $p$, such that every smooth solution $v$ of
\begin{equation*}
 \diver((\tau+|Dv|^2)^{(p-2)/2}Dv)=0,\qquad \tau>0,
\end{equation*}
on a neighborhood of $\overline{B_{2R}}\subset\mathbb R^3$ satisfies
\begin{equation}\label{eq:jin-energy}
 \int_{B_R\cap\{|Dv|>0\}}\frac{|D^2v|^2}{|Dv|}\,dx
 \leq L_p R^{-2}\int_{B_{2R}}|Dv|\,dx .
\end{equation}
The constant does not depend on $\tau$ or $v$. The denominator in this estimate is the actual gradient magnitude $|Dv|$, not the regularized magnitude $(\tau+|Dv|^2)^{1/2}$. It is precisely this weighted estimate, rather than just its unweighted $W^{2,2}$ consequence, that permits the following transfer.

\begin{proposition}[Conditional square-root gradient transfer in dimension three]\label{prop:jin-transfer}
Let $n=3$ and $5\leq p<\infty$, and assume \eqref{eq:jin-energy} for this $p$, with $L_p$ independent of the regularization parameter and the solution. Every scalar $p$-harmonic function has $N_1\in W^{1,2+\varepsilon_{\rm hi}}_{\mathrm{loc}}(\Omega)$ for some $\varepsilon_{\rm hi}=\varepsilon_{\rm hi}(p)>0$. In particular,
\begin{equation*}
 \int_{B_r}|DN_1|^2\,dx
 \leq Cr^{-2}\int_{B_{4r}}|N_1-(N_1)_{B_{4r}}|^2\,dx
 \qquad (B_{4r}\Subset\Omega).
\end{equation*}
\end{proposition}

\begin{proof}
Fix a concentric pair $B_r\subset B_{2r}\Subset\Omega$ and an intermediate slightly larger ball $B_{2r+\rho_{\rm ext}}\Subset\Omega$. By the standard local regularized Dirichlet approximation for scalar $p$-harmonic functions (see \cite[Section 3]{Sarsa}), there are smooth solutions $v_k$ on $B_{2r+\rho_{\rm ext}}$ with the boundary values of $u$ and regularization parameters $\tau_k\downarrow0$, such that $Dv_k\to Du$ uniformly on $\overline{B_{2r}}$. This approximation uses only first-order regularity.

Write $Q_k=|Dv_k|$. For a second, independent parameter $d>0$ define
\begin{equation*}
 Z_{k,d}=(d^2+Q_k^2)^{-1/4}Dv_k.
\end{equation*}
For the map $z\mapsto (d^2+|z|^2)^{-1/4}z$, the radial and tangential eigenvalues of the derivative are $(d^2+|z|^2)^{-1/4}(1-|z|^2/(2(d^2+|z|^2)))$ and $(d^2+|z|^2)^{-1/4}$, respectively. Hence, on $\{Q_k>0\}$,
\begin{equation*}
 |DZ_{k,d}|^2\leq (d^2+Q_k^2)^{-1/2}|D^2v_k|^2
 \leq Q_k^{-1}|D^2v_k|^2.
\end{equation*}
Also $D^2v_k=0$ almost everywhere on $\{Dv_k=0\}$ by Sobolev locality of the gradients of the smooth functions $v_k$. Thus there is no additional integral over that set. Applying \eqref{eq:jin-energy} and using $|Z_{k,d}|^2\leq Q_k$ gives, uniformly in $d$,
\begin{equation*}
 \int_{B_r}|DZ_{k,d}|^2\,dx
 \leq L_p r^{-2}\int_{B_{2r}}Q_k\,dx.
\end{equation*}
For fixed $k$, as $d\downarrow0$ the fields $Z_{k,d}$ converge strongly in $L^2(B_r)$ to $Q_k^{-1/2}Dv_k$ with value zero at $Q_k=0$, by pointwise convergence and domination of the squared norm by $Q_k$. Weak compactness and lower semicontinuity give the same estimate for the weak derivatives of this limit field. Next use the locally uniform convergence of $Dv_k$ and the continuity of the square-root radial map at zero. The fields $Q_k^{-1/2}Dv_k$ converge uniformly on $\overline{B_r}$ to $N_1$, and $Q_k\to Q$ uniformly on $\overline{B_{2r}}$. Applying weak compactness once more yields
\begin{equation*}
 N_1\in W^{1,2}(B_r),\qquad
 \int_{B_r}|DN_1|^2\,dx\leq L_p r^{-2}\int_{B_{2r}}Q\,dx .
\end{equation*}
Since the pair of balls was arbitrary, this is exactly \eqref{eq:unanchored-input} at $s=1$. Theorem \ref{thm:anchored-native} gives the stated anchored estimate and higher integrability.
\end{proof}

\begin{corollary}[Conditional extension to all finite exponents in dimension three]\label{cor:three-all}
Assume \eqref{eq:jin-energy} for every finite $p\geq5$, with $L_p$ depending only on $p$. Then, for every finite $p>2$, each scalar $p$-harmonic function in an open subset of $\mathbb R^3$ belongs to $W^{2,2+\varepsilon_{\rm hi}}_{\rm loc}$ for some $\varepsilon_{\rm hi}=\varepsilon_{\rm hi}(p)>0$.
\end{corollary}

\begin{proof}
For $2<p<5$ this follows independently of Jin from Corollary \ref{cor:unweighted-higher}. For $5\leq p<\infty$, Proposition \ref{prop:jin-transfer} gives $N_1\in W^{1,2+\varepsilon_{\rm hi}}_{\rm loc}$. The original gradient is $Du=|N_1|N_1$, a locally Lipschitz composition of the locally bounded field $N_1$. The Sobolev chain rule gives the assertion.
\end{proof}

Independently of the external hypothesis, Corollaries~\ref{cor:unweighted-higher} and \ref{eff:delta} already give the unconditional range $2<p\leq5+\delta_3$ in dimension three. The additional hypothesis extends the conclusion beyond this range. Its curvature argument is not part of the present proof. No claim of optimality, novelty, or priority is made for the additional three-dimensional conclusion, and no extension of that conclusion to bounded nonzero sources is asserted.

\section{A quantitative choice of the exponent}\label{sec:effective-exponent}

The preceding argument gives a strict gain below $m=(p-1)(n-2)/(n-1)$ by compactness. We show that the gain can be chosen by finite expressions involving only $n$ and $p$. The constants obtained below are not intended to be sharp. We replace the charge in Lemma~\ref{lem:charge} by a charge on a larger, explicitly specified interior ball. This permits a direct use of the differential equation of the unit normal field, without quantitative analytic continuation.

Let $P=p-1$. For a compact interval $I=[s_-,s_+]$ with $0<s_-<m<s_+<2P$, we shall specify positive numbers $c_{\rm hom},K_{\rm loc},C_{\rm ov}$. The exponent will be
\begin{equation}\label{eff:main-exp}
 \sigma=\min\left\{\frac m4,\frac{3s_-c_{\rm hom}}{128P K_{\rm loc}C_{\rm ov}}\right\},
 \qquad s_*=m-\sigma .
\end{equation}
For $I=[m/2,3m/2]$, all the quantities in this formula are explicitly determined by $n,p$. The same choice is valid in the bounded-source theorem; the smallness threshold determining the natural source scale still arises from qualitative compactness.

\subsection{An effective first-order estimate}

We first record a conservative quantitative form of the interior gradient estimate for the homogeneous equation. The argument uses elementary De Giorgi iteration. Put $\omega_n=|B_1|$, $\gamma=n/(n-2)$, and define
\begin{align}
 K&=P4^{p-2},\quad S_n=16n^2,\quad C_{\rm Poin}=2^{n+1},\quad B_d=16^\gamma,\notag\\
 A_d&=\bigl(2^{20}S_n K\max\{1,\omega_n^{2/n}\}\bigr)^\gamma,\notag\\
 \nu_d&=\frac12\min\{1,A_d^{-1/(\gamma-1)}B_d^{-1/(\gamma-1)^2}\},
 \quad \rho_d=(1-\nu_d/2)^{1/n},\notag\\
 C_d&=\frac{16 C_{\rm Poin}\sqrt K}{\nu_d(1-\rho_d)},\quad
 N_d=\left\lceil\left(\frac{2C_d}{\nu_d\rho_d^n}\right)^2\right\rceil+1,\notag\\
 d_d&=2^{-N_d-2},\quad \alpha=\frac{-\log(1-d_d)}{\log4}.\label{eff:alpha}
\end{align}
In particular, $0<\alpha<1$. We use the Sobolev inequality
\begin{equation*}
 \|w\|_{L^{2\gamma}(\R^n)}^2\le S_n\int_{\R^n}|Dw|^2\,dx
\end{equation*}
for compactly supported Sobolev functions. The constant $S_n$ is more than sufficient: writing $|\varphi|\le\int|\partial_i\varphi|\,dx_i$ for each coordinate and applying iterated H\"older inequality gives $\|\varphi\|_{n/(n-1)}\le\prod_i\|\partial_i\varphi\|_1^{1/n}\le\|D\varphi\|_1$. Apply this to $\varphi=|w|^{2(n-1)/(n-2)}$ and use Cauchy--Schwarz. We also use
\begin{equation*}
 \int_{B_r}\int_{B_r}|w(x)-w(y)|\,dx\,dy
 \le C_{\rm Poin} r |B_r|\int_{B_r}|Dw|\,dx .
\end{equation*}
Indeed, integrate the line-segment identity for $w(x)-w(y)$ over pairs, splitting the segment parameter at $1/2$ and changing the variable $x$ on one half and $y$ on the other. On each half the Jacobian cost is at most $2^n$.

\begin{lemma}[Effective homogeneous gradient modulus]\label{eff:holder-lemma}
For a scalar $p$-harmonic function in $B_R(z)$,
\begin{equation}\label{eff:centred-holder}
 |Du(y)-Du(z)|\le64\sup_{B_R(z)}Q\,(|y-z|/R)^\alpha,\qquad y\in B_R(z).
\end{equation}
In particular, for $x,y\in B_{R/2}(z)$,
\begin{equation}\label{eff:pair-holder}
 |Du(x)-Du(y)|\le128\sup_{B_R(z)}Q\,(|x-y|/R)^\alpha.
\end{equation}
For $p$ in a compact subinterval of $(2,\infty)$, the same estimates hold with $\alpha$ calculated at the upper endpoint.
\end{lemma}

\begin{proof}
We begin with the smooth regularized homogeneous equation on $B_1$,
\begin{equation*}
 \diver\big((\tau+|Du|^2)^{(p-2)/2}Du\big)=0,\qquad \tau>0.
\end{equation*}
For a unit vector $\omega$ put $M=(\tau+\sup_{B_1}|Du|^2)^{1/2}$, $v=\omega\cdot Du/M$. Then $|v|\le1$ and $\diver(A(x)Dv)=0$, where
\begin{equation*}
 A(x)=M^{-(p-2)}D_\xi\big((\tau+|\xi|^2)^{(p-2)/2}\xi\big)\big|_{\xi=Du(x)} .
\end{equation*}
The matrix is positive definite with $A\le P I$, and $A\ge4^{-(p-2)}I$ on the set $\{v\ge1/4\}$. Taking $\eta^2(v-k)_+$, $k\ge1/4$, as a test function gives
\begin{equation}\label{eff:DG-up}
 \int\eta^2|D(v-k)_+|^2\,dx\le4K\int(v-k)_+^2|D\eta|^2\,dx .
\end{equation}
For $1/4\le l<k\le1/2$ put $w=\min\{(k-v)_+,k-l\}$. Taking the full negative-part test $\eta^2(k-v)_+$, we obtain
\begin{equation}\label{eff:DG-down}
 \int\eta^2|Dw|^2\,dx\le16K\int_{\{v<k\}}|D\eta|^2\,dx .
\end{equation}
Indeed, the untruncated energy in the $A$-metric is at most $4\int(k-v)_+^2 A D\eta\cdot D\eta$; $Dw=0$ outside $\{l<v<k\}$, on which the lower ellipticity bound is valid. No such lower bound is needed on $\{v\le l\}$.

Let $r_j=1/2+2^{-j-1}$ and take cutoffs equal to one in $B_{r_{j+1}}$, supported in $B_{r_j}$, with $|D\eta_j|\le2^{j+3}$. Set $k_j=1/4+2^{-j-2}$, $E_j=B_{r_j}\cap\{v<k_j\}$ and $w_j=\min\{(k_j-v)_+,k_j-k_{j+1}\}$. By \eqref{eff:DG-down}, the product rule and the Sobolev inequality,
\begin{equation*}
 2^{-2j-6}|E_{j+1}|^{1/\gamma}\le
 S_n\int|D(\eta_jw_j)|^2\,dx\le2^{12}S_nK4^j |E_j|.
\end{equation*}
Thus $Y_j=|E_j|/\omega_n$ satisfies $Y_{j+1}\le A_d B_d^jY_j^\gamma$. The same recurrence holds for the upper-level sets $E_j=B_{r_j}\cap\{v>k_j\}$, with $k_j=k+(a/2)(1-2^{-j})$, $a=1-k$, $1/2\le k<1$. In this case use \eqref{eff:DG-up}, observe that $(v-k_j)_+\le a$, and use the level gap $a2^{-j-2}$ in the Sobolev inequality; the factor $a^2$ cancels. If $Y_0\le\nu_d$, substitution into the recurrence gives
\begin{equation*}
 Y_j\le A_d^{-1/(\gamma-1)}B_d^{-1/(\gamma-1)^2}
           B_d^{-j/(\gamma-1)}2^{-\gamma^j},
\end{equation*}
which tends to zero. We have proved
\begin{align}
 |\{v<1/2\}\cap B_1|\le\nu_d\omega_n&\quad\Longrightarrow\quad v\ge1/4
       \text{ in }B_{1/2},\label{eff:DG-L}\\
 |\{v>k\}\cap B_1|\le\nu_d\omega_n&\quad\Longrightarrow\quad v\le(1+k)/2
       \text{ in }B_{1/2}.\label{eff:DG-U}
\end{align}
The proofs are invariant under spatial rescaling.

Suppose instead that $|\{v\le1/2\}\cap B_1|\ge\nu_d\omega_n$. The low set in $B_{\rho_d}$ has measure at least $(\nu_d/2)\omega_n$. Take $k_j=1-2^{-j-1}$, $a_j=k_{j+1}-k_j=2^{-j-2}$ and $A_j=B_{\rho_d}\cap\{v>k_j\}$. Applying the pairwise Poincar\'e inequality to $\min\{(v-k_j)_+,a_j\}$ gives
\begin{equation*}
 a_j|A_{j+1}||\{v\le1/2\}\cap B_{\rho_d}|
 \le C_{\rm Poin}\rho_d |B_{\rho_d}|
        \int_{B_{\rho_d}\cap\{k_j<v<k_{j+1}\}}|Dv|\,dx .
\end{equation*}
By \eqref{eff:DG-up} with cutoff gradient at most $2/(1-\rho_d)$, and by $1-k_j=2a_j$,
\begin{equation*}
 \int_{B_{\rho_d}\cap\{k_j<v<k_{j+1}\}} |Dv|^2\,dx
 \le\frac{64K a_j^2\omega_n}{(1-\rho_d)^2}.
\end{equation*}
Cauchy--Schwarz yields $|A_{j+1}|/\omega_n\le C_d(|A_j\setminus A_{j+1}|/\omega_n)^{1/2}$. Level equality sets cause no difficulty by Sobolev locality. Summing the squares over $0\le j<N_d$ and using nesting gives $|A_{N_d}|/\omega_n\le C_d/\sqrt{N_d}\le\nu_d\rho_d^n/2$. Apply \eqref{eff:DG-U} on $B_{\rho_d}$ with level $k_{N_d}$. Since $\rho_d/2>1/4$, we conclude $v\le1-d_d$ in $B_{1/4}$.

The estimates apply to every direction, with uniform constants. Passing to the actual solution by the usual regularized Dirichlet approximation on a slightly larger interior ball gives, for $M_R=\sup_{B_R(z)}Q$,
\begin{equation}\label{eff:dichotomy}
 \begin{cases}
 Q\ge M_R/4\quad\text{in }B_{R/2}(z), &\text{or}\\
 \sup_{B_{R/4}(z)}Q\le(1-d_d)M_R .
 \end{cases}
\end{equation}
Here the approximating gradients converge locally uniformly by classical first-order regularity and Dirichlet comparison (as in \cite[Section 3]{Sarsa}); no numerical estimate from that convergence is used. If some directions satisfy \eqref{eff:DG-L} along a subsequence, pass to a convergent subsequence of directions. Otherwise the bound $v\le1-d_d$ for every direction gives the second alternative by taking the direction of the gradient itself. Both assertions pass to the limit. Increasing interior radii gives the same dichotomy for the ball in the statement. The cases $M_R=0,\infty$ are immediate or vacuous.

On a ball in the first alternative, every directional derivative solves a uniformly elliptic divergence equation on the half-radius ball, with ellipticity ratio at most $K$. For a derivative normalized by its infimum and oscillation, either it or its complement has low-set measure at least half the ball. The upper-level argument above, now valid at every level by uniform ellipticity, contracts the oscillation by $1-d_d$ at quarter scales.

Iterate \eqref{eff:dichotomy} at radii $R4^{-j}$ about $z$. Until the first ball in the first alternative, the gradient supremum is at most $M_R4^{-\alpha j}$. If the first alternative occurs first at $j_0$, directional derivative oscillations from that half-radius ball onwards decay by the same power under quarter-scale iteration. For $d=|y-z|\le R/8$, choose $j$ so that $R4^{-j}/8\ge d>R4^{-j-1}/8$. If no first alternative has occurred at scales less than $j$, both gradient values have magnitude at most $M_R(32d/R)^\alpha$. Otherwise the directional derivative difference is at most $2M_R(8d/R)^\alpha$. Taking the supremum over directions and using $\alpha\le1$ proves \eqref{eff:centred-holder}; the distances $d>R/8$ follow from the trivial bound $2M_R$. For \eqref{eff:pair-holder}, apply the centered estimate on $B_{R/2}(x)$ when $|x-y|<R/2$ and the trivial bound otherwise. Uniformity follows by using $K=(p_+-1)4^{p_+-2}$ on an interval with upper endpoint $p_+$.
\end{proof}

\subsection{Effective derivative bounds on positive-gradient balls}

Only smooth solutions of the homogeneous equation on balls with positive lower gradient magnitude are used in this subsection. We first specify the constants in an elementary linear estimate by comparison with the frozen-coefficient equation.

\begin{lemma}[Effective linear gradient estimate]\label{eff:linear}
Suppose $v,\mathsf A,G$ are smooth on $B_1$ and
\begin{equation*}
 \diver(\mathsf A(x)Dv)=\diver G,\quad
 0<\lambda_{\rm ell} I\le \mathsf A\le\Lambda_{\rm ell} I,\quad
 [\mathsf A]_a\le L_{\mathsf A},\quad [G]_a\le L_G,\quad \|v\|_\infty\le M ,
\end{equation*}
where $\mathsf A$ is symmetric and $0<a\le1$. Define $K_0=\Lambda_{\rm ell}/\lambda_{\rm ell}$, $\beta=a/2$, and
\begin{align}
 C_f&=4(n+2)4^{n+2}K_0^n,\quad
 \theta=\min\{1/8,1/(4\sqrt{K_0}),(4C_f)^{-1/(1-\beta)}\},\notag\\
 C_1^f&=C_f+\theta^{-n/2},\quad
 C_m=1+\theta^{-n/2}/(1-\theta^\beta),\notag\\
 r_f&=\min\left\{\frac18,\left(\frac{\lambda_{\rm ell}\theta^\beta}
                      {4C_1^f(1+C_m)(1+L_{\mathsf A})}\right)^{1/a}\right\},\notag\\
 J_f&=2^{n/2+6}\left(\frac{K_0M}{r_f}+\frac{(2r_f)^a L_G}{\lambda_{\rm ell}}\right),\quad
 D_f=2J_f+4C_1^f\theta^{-\beta}\lambda_{\rm ell}^{-1}r_f^a L_G.\label{eff:lin-const}
\end{align}
Then
\begin{align}
 \|Dv\|_{C^{0,\beta}(B_{1/2})}
 &\le{\mathfrak L}(n,a,\lambda_{\rm ell},\Lambda_{\rm ell},L_{\mathsf A},L_G,M)\notag\\
 &:=D_f\left(C_m+(2C_m+2^{1+n/2})(2/(\theta r_f))^\beta\right).
 \label{eff:lin-est}
\end{align}
The H\"older norm is the sum of the supremum and seminorm. For the corresponding equation on a ball $B_R$, the right-hand side after rescaling is bounded by $2\max\{R^{-1},R^{-1-\beta}\}{\mathfrak L}$ with arguments $R^aL_{\mathsf A},R^{1+a}L_G,M$ in place of $L_{\mathsf A},L_G,M$.
\end{lemma}

\begin{proof}
For $z\in B_{1/2}$ take a cutoff on $B_{2r_f}(z)$, equal to one on $B_{r_f}(z)$, with gradient at most $2/r_f$. Test by $\eta^2v$ and replace $G$ by $G-G(z)$. With $X=\eta|Dv|$, $Y=|v||D\eta|$ and $W=\lambda_{\rm ell}^{-1}\eta|G-G(z)|$, one obtains
\begin{equation*}
 \int X^2\le2K_0\int XY+\int XW+2\int WY
 \le\frac12\int X^2+(4K_0^2+1)\int Y^2+2\int W^2.
\end{equation*}
Thus $(\avg_{B_{r_f}(z)}|Dv|^2)^{1/2}\le J_f$.

On any $B_r(z)\subset B_1$, let $w$ solve $\diver(\mathsf A(z)Dw)=0$ with trace $v$. Test the difference equation by $v-w$ and subtract the constant $G(z)$, obtaining
\begin{equation}\label{eff:frozen-diff}
 \left(\avg_{B_r}|D(v-w)|^2\right)^{1/2}
 \le\lambda_{\rm ell}^{-1}r^a\left(L_{\mathsf A}\left(\avg_{B_r}|Dv|^2\right)^{1/2}+L_G\right).
\end{equation}
For every constant vector $c$ the frozen equation gives
\begin{equation}\label{eff:frozen-osc}
 \left(\avg_{B_{\theta r}}|Dw-(Dw)_{B_{\theta r}}|^2\right)^{1/2}
 \le C_f\theta\left(\avg_{B_r}|Dw-c|^2\right)^{1/2}.
\end{equation}
To see this, apply $x=z+(\mathsf A(z)/\Lambda_{\rm ell})^{1/2}y$. The transformed domain contains $B_r$ and the points of $B_{\theta r}(z)$ satisfy $|y|\le\sqrt{K_0}\theta r\le r/4$. Each component of $Dw-c$ in these variables is harmonic. Use the harmonic mean value identity with the radial probability kernel $\psi(y)=(n+2)(1-|y|^2)_+/(2\omega_n)$ at scale $r/4$, then differentiate the convolution. Its gradient has $L^2$ norm at most $(n+2)\omega_n^{-1/2}(r/4)^{-1-n/2}$. Cauchy--Schwarz and the change of variables bound $\sup_{B_{\theta r}}|D^2w|$ by
\begin{equation*}
 (n+2)4^{1+n/2}K_0^{1/2+n/4}r^{-1}
                    \left(\avg_{B_r}|Dw-c|^2\right)^{1/2}.
\end{equation*}
Integrating the derivative along segments in the smaller ball proves \eqref{eff:frozen-osc}, with room in the choice of $C_f$.

Set $r_j=r_f\theta^j$, $E_j=(\avg_{B_{r_j}(z)}|Dv-(Dv)_{B_{r_j}(z)}|^2)^{1/2}$ and $m_j=|(Dv)_{B_{r_j}(z)}|$. Combining \eqref{eff:frozen-diff}, \eqref{eff:frozen-osc}, and the minimizing property of the average, we obtain
\begin{align*}
 E_{j+1}&\le C_f\theta E_j+C_1^f\lambda_{\rm ell}^{-1}r_j^a[L_{\mathsf A}(E_j+m_j)+L_G],\\
 m_{j+1}&\le m_j+\theta^{-n/2}E_j.
\end{align*}
Observe that $E_0+m_0\le D_f$, $C_f\theta\le\theta^\beta/4$, and $C_1^f\lambda_{\rm ell}^{-1}r_f^a L_{\mathsf A}\le\theta^\beta/[4(1+C_m)]$. It follows inductively that $E_j\le D_f\theta^{\beta j}$ and $m_j\le C_mD_f$. Indeed, the latter bound follows by summing the mean-difference inequality, and the three terms in the excess inequality are each bounded by $D_f\theta^{\beta(j+1)}/4$, using $a\ge\beta$ and the definition of $D_f$. Summing the remaining mean differences also gives $|Dv(z)-(Dv)_{B_{r_j}(z)}|\le C_mD_f\theta^{\beta j}$.

For $z,z'\in B_{1/2}$ at distance $d\le r_f/2$, choose $j$ with $r_j\ge2d>\theta r_j$. The two averaging balls have in common a ball of radius $r_j/2$. Comparing their averages on this common ball by Cauchy--Schwarz bounds their difference by $2^{1+n/2}D_f\theta^{\beta j}$. Compare the two averages to their point values and use $r_j<2d/\theta$, obtaining
\begin{equation*}
 |Dv(z)-Dv(z')|\le(2C_m+2^{1+n/2})D_f(2d/(\theta r_f))^\beta.
\end{equation*}
For larger distances the supremum bound suffices. This proves \eqref{eff:lin-est}. Under the spatial rescaling $x=x_0+Ry$, the right-hand side vector becomes $RG(x_0+Ry)$ and the derivative rescales by $R^{-1}$, giving the last assertion.
\end{proof}

\begin{lemma}[Effective positive-gradient derivative bounds]\label{eff:warm}
Let $U$ be scalar $p$-harmonic on $B_1$ with $q_0\le|DU|\le A$, $0<q_0<1\le A$, and $[DU]_a\le H_1$, $0<a\le1$. Put
\begin{align}
 \lambda_{\rm ell}&=q_0^{p-2},\quad \Lambda_{\rm ell}=P A^{p-2},\quad
 M_c=10^8 n^{10}(p+20)^{20}(1+A+2/q_0)^{p+20},\notag\\
 C_c&=(1+4/q_0)M_c,\quad C_0=10n^4 C_c(1+2H_1),\quad a_i=a/2^i,\notag\\
 J_2&=1+2n^2{\mathfrak L}(n,a,\lambda_{\rm ell},\Lambda_{\rm ell},C_0,0,A),\notag\\
 J_3&=1+2n^3\max\{2,2^{1+a_2}\}
          {\mathfrak L}(n,a_1,\lambda_{\rm ell},\Lambda_{\rm ell},C_0,C_0J_2^2,J_2),\notag\\
 J_4&=1+2n^4\max\{4,4^{1+a_3}\}
          {\mathfrak L}(n,a_2,\lambda_{\rm ell},\Lambda_{\rm ell},C_0,
                       C_0(1+J_2+J_3)^3,J_3),\notag\\
 K_w&=100n^4(1+C_c)^3(1+J_2+J_3+J_4)^4.\label{eff:warm-def}
\end{align}
Then $|Dg|,|D^2g|\le K_w$ on $B_{1/8}$, where $g=D_e e$ and $e=DU/|DU|$. On a ball of radius $R_w$, with these hypotheses after the unit spatial rescaling preserving the gradient amplitude, the bounds are $|Dg|\le R_w^{-2}K_w$ and $|D^2g|\le R_w^{-3}K_w$.
\end{lemma}

\begin{proof}
On the positive-gradient ball classical nondegenerate regularity gives a smooth solution. Write $r=|\xi|$ and $T=1+A+2/q_0$. The differentiated matrix and the unit direction are
\begin{equation*}
 \mathcal{A}(\xi)=r^{p-2}I+(p-2)r^{p-4}\xi\otimes\xi,
 \qquad e(\xi)=r^{-1}\xi.
\end{equation*}
We verify explicitly that $M_c$ bounds their derivative tensors through order four on $q_0/2\leq r\leq A+1$. For $b\in\{p-2,p-4,-1\}$ and $0\leq k\leq4$, an indexed derivative of $r^b$ is a sum over partitions of the $k$ differentiation indices into singletons and pairs. The numbers of such partitions are $1,1,2,4,10$ for $k=0,1,2,3,4$. If there are $h$ blocks, the coefficient is a product of $h$ factors from $b,b-2,b-4,b-6$. Its absolute value is at most $(p+10)^4$, and the coordinate factors leave the radial power $r^{b-k}$. Since $r,r^{-1}\leq T$ and $|b-k|\leq p+8$ on this annulus, every indexed derivative satisfies
\begin{equation*}
 |\partial_{i_1}\cdots\partial_{i_k}r^b|
 \leq10(p+10)^4T^{p+8}.
\end{equation*}
For each component of $r^{p-4}\xi\otimes\xi$, Leibniz's rule has at most $3^k\leq81$ assignments of derivatives to the three factors, and the two coordinate factors and their derivatives contribute at most $T^2$. Thus each component of $D_\xi^k\mathcal{A}$ is bounded by $820(p+10)^5T^{p+10}$. The corresponding bound for each component of $D_\xi^ke$ is $160(p+10)^4T^{p+9}$, since there are at most $2^k\leq16$ assignments. A matrix derivative of order at most four has at most $n^6$ components, and a vector derivative has at most $n^5$. Hence their Euclidean tensor norms obey
\begin{equation*}
 \|D_\xi^k\mathcal{A}\|+\|D_\xi^ke\|
 \leq1000n^3(p+10)^5T^{p+10}\leq M_c,
 \qquad 0\leq k\leq4.
\end{equation*}

For vectors of mutual distance at most $q_0/2$ in the original annulus the segment stays in the padded annulus, permitting use of the derivative bound; otherwise use twice the supremum divided by $q_0/2$. Consequently the derivatives of order at most three have supremum and Lipschitz bounds $C_c$ on the original annulus. In particular, the $C^{0,a'}$ norms of $(D_\xi^j\mathcal{A})(DU)$ are bounded by $C_c(1+2H_1)$ for $j\le3$ and $0<a'\le a$, and likewise for the derivatives of $e(\xi)$.

Let $\mathsf A_U=\mathcal{A}(DU)$. Each $U_i$ satisfies the homogeneous linear divergence equation with coefficient $\mathsf A_U$. Applying Lemma~\ref{eff:linear} gives $\|D^2U\|_{C^{0,a_1}(B_{1/2})}\le J_2$. Differentiating once more gives
\begin{equation*}
 \diver(\mathsf A_U DU_{ij})=-\diver\big((\partial_j\mathsf A_U)DU_i\big).
\end{equation*}
On $B_{1/2}$, the right-hand side vector has $C^{0,a_1}$ norm bounded by $C_0J_2^2$, by the chain rule for $\partial_j\mathsf A_U$ and the product inequality for the H\"older norm. Use the rescaled estimate of Lemma~\ref{eff:linear} to obtain $\|D^3U\|_{C^{0,a_2}(B_{1/4})}\le J_3$. The arguments of $\mathfrak L$ in \eqref{eff:warm-def} overestimate the coefficient and right-hand side seminorms after rescaling.

Similarly,
\begin{equation*}
 \diver(\mathsf A_U DU_{ijk})=-\diver\big((\partial_{jk}\mathsf A_U)DU_i
       +(\partial_j\mathsf A_U)DU_{ik}+(\partial_k\mathsf A_U)DU_{ij}\big).
\end{equation*}
Since $D^2\mathsf A_U$ is the sum of contractions of $D_\xi^2\mathcal{A}(DU)$ with two Hessians, and of $D_\xi\mathcal{A}(DU)$ with the third derivative, the right-hand side vector is bounded in $C^{0,a_2}$ by $C_0(1+J_2+J_3)^3$. The linear estimate on $B_{1/4}$ gives $\|D^4U\|_{C^{0,a_3}(B_{1/8})}\le J_4$. The factors $n^2,n^3,n^4$ and $10n^4$ in \eqref{eff:warm-def} dominate the finite sums of scalar component estimates and tensor contractions.

The chain rule for $e=e(DU)$ now yields on $B_{1/8}$
\begin{align*}
 |De|&\le C_cJ_2,\qquad |D^2e|\le2C_c(J_2+J_3)^2,\\
 |D^3e|&\le8C_c(1+J_2+J_3+J_4)^3 .
\end{align*}
As $g=D_e e$, one has $|Dg|\le |D^2e|+|De|^2$ and $|D^2g|\le|D^3e|+3|D^2e||De|$. These are bounded by $K_w$. The scaling follows by the chain rule under $U(y)=(u(x_0+R_wy)-u(x_0))/R_w$.
\end{proof}

\subsection{A direct quantitative geometric charge}

Take $H_0=256$, twice the bound in \eqref{eff:pair-holder}. Define
\begin{align}
 b_*&=\tfrac14\log(8/7),\quad \nu_*=\frac{P^2b_*^2}{2(n-1)},\quad
 \ell_*=3+\frac{4(n-1)}{Pb_*},\quad A_*=(\ell_*+1)^{2/s_-},\notag\\
 r_0&=\min\left\{\frac1{64},\frac18[16H_0(A_*+1)]^{-1/\alpha}\right\},\quad
 R_w=8r_0.\label{eff:geometry-const}
\end{align}
In \eqref{eff:warm-def} choose $q_0=3/4$, $A=A_*$, $a=\alpha$ and $H_1=H_0A_*$. Then put
\begin{align}
 K_1&=1+R_w^{-2}K_w,\qquad K_2=1+R_w^{-3}K_w,\notag\\
 t_0&=\min\{r_0/8,\nu_* /(2nK_2)\},\quad
 \eta_*=\min\{b_*,\nu_* t_0/(4n)\},\quad
 d_*=\min\{r_0/8,\eta_* /(2K_1)\},\notag\\
 c_{\rm hom}&=(3/4)^{s_+}\omega_n d_*^n\eta_*^2/4.\label{eff:charge-const}
\end{align}

\begin{lemma}[Effective enlarged-charge lemma]\label{eff:charge}
Let $u$ be scalar $p$-harmonic in $B_{\ell_*+1}$ with $Q(0)\le1/2$, $\max_{\overline B_1}Q=1$ and $\sup_{B_{\ell_*+1}}Q\le A_*$. For each $s\in I$,
\begin{equation}\label{eff:charge-est}
 \int_{B_{\ell_*}\cap\{Q>3/4\}}Q^s|g|^2\,dx\ge c_{\rm hom}.
\end{equation}
\end{lemma}

\begin{proof}
At every $x\in B_{\ell_*-1}$ with $Q(x)\ge7/8$, apply \eqref{eff:pair-holder} on $B_1(x)$. By \eqref{eff:geometry-const}, one has $13/16\le Q\le A_*$ on $B_{R_w}(x)$. After rescaling to the unit ball while preserving the gradient amplitude, the normalized gradient seminorm is at most $H_0A_*$. Lemma~\ref{eff:warm} gives $|Dg|\le K_1$, $|D^2g|\le K_2$ on $B_{r_0}(x)$.

Suppose \eqref{eff:charge-est} fails. Then $|g|\le\eta_*$ on $B_{r_0/4}(x)$. Otherwise the bound for $Dg$ yields $|g|>\eta_*/2$ on a ball of radius $d_*$ centered at a violating point, contained in $B_{r_0/2}(x)\subset\{Q>3/4\}$. Its contribution to \eqref{eff:charge-est} is at least $c_{\rm hom}$. Taylor's formula with increment $t_0$ in each coordinate direction gives
\begin{equation*}
 |\partial_i g_i(x)|\le2\eta_*/t_0+K_2t_0/2\le\nu_*/n.
\end{equation*}
Thus $|g|\le\eta_*$ and $|\diver g|\le\nu_*$ on $B_{\ell_*-1}\cap\{Q\ge7/8\}$.

Recall that $h=D_e\log Q$ and $\tr S=-Ph$ in the homogeneous equation. Since $D_e e=g$ and $\diver e=\tr S$, differentiation gives
\begin{equation*}
 -P D_e h=D_e\diver e
 =\diver(D_e e)-\tr((De)^2)=\diver g-|S|^2.
\end{equation*}
Here the adapted-frame matrix of $De$ has normal row zero, normal column $(0,g)^T$ and tangential block $S$, so $\tr((De)^2)=|S|^2$. Using $|S|^2\ge P^2h^2/(n-1)$, we obtain
\begin{equation}\label{eff:riccati}
 D_e h\ge\frac{P}{n-1}h^2-\frac{\nu_*}{P}
\end{equation}
on the high-gradient region.

If $|h(x)|>b_*$ at a point $x\in\overline B_1$, $Q(x)\ge7/8$, consider the integral curve of $\epsilon e$ with $\epsilon=\operatorname{sign}h(x)$ fixed. Along the curve let $k=\epsilon h$. Then $k'=D_e h$ and $(\log Q)'=k$. While $k>b_*$, \eqref{eff:riccati} implies $k'\ge Pk^2/(2(n-1))$. Thus $k$ increases and $Q$ increases; the curve cannot leave the high-gradient region through a lower speed. Integration of the inequality for $1/k$ forces blow-up by time $2(n-1)/(Pb_*)<\ell_*-2$. Since the curve has unit speed and starts in $\overline B_1$, it stays in a compact interior part of $B_{\ell_*-1}$ during this time. Its positive lower speed keeps it in a compact positive-gradient region where $h$ is bounded and the curve continues. This contradiction proves $|h|\le b_*$ on $\overline B_1\cap\{Q\ge7/8\}$.

Join a maximum point in $\overline B_1$ to the origin by a segment. Before its first drop from $1$ to $7/8$ one has $|D\log Q|\le|h|+|g|\le2b_*$. The length is at most one, so the logarithmic drop would be at most $2b_*=\frac12\log(8/7)$. This contradicts the first drop; such a drop must occur because $Q(0)\le1/2$.
\end{proof}

\subsection{The numerical covering coefficient and absorption}

We modify the proof of Lemma~\ref{lem:korn} by using the enlarged charge. Define
\begin{align}
 \kappa_*&=\min\{1/40,(30H_0)^{-1/\alpha}/20\},\quad
 c_*=\min\{1/2,[3/(8H_0(A_*+1))]^{1/\alpha}\},\notag\\
 C_*&=\max\{2(\ell_*+1)+1,[4H_0(2A_*+1)]^{1/\alpha}\},\qquad
 L_*=2+C_*+\ell_*,\notag\\
 \mu&=(9/10)^{p-2},\quad B_* = P(11/10)^{p-2},\quad
 W_*=\max\{1,(9/10)^{s_--2},(11/10)^{s_+-2}\},\notag\\
 C_H&=2\omega_n(10\kappa_*)^n\left[
 \frac{4B_*^2(121/100)/\mu^2+1/\mu}{25\kappa_*^2}
          +\frac n{\mu^2}+\frac{121/100}{\mu}\right],\quad
 K_{\rm loc}=W_*C_H,\notag\\
 C_{\rm ov}&=\left[\frac{C_*(\ell_*+\kappa_*)}{c_*\kappa_*}\right]^n,\quad
 C_{\rm cov}=\frac{2K_{\rm loc}C_{\rm ov}}{c_{\rm hom}},\quad
 C_{\rm add}=\frac{16 L_*^2 K_{\rm loc}}{\omega_n\kappa_*^n}.
 \label{eff:overlap-const}
\end{align}

\begin{lemma}[Quantitative weighted covering]\label{eff:cover}
Let $B_r\subset B_t\subset B_R\Subset\Omega$ be concentric, with $M=\sup_{B_R}Q$ and $\chi_\varepsilon$ as in Subsection~\ref{sec:covering}. For a homogeneous solution and every $s\in I$,
\begin{equation}\label{eff:cover-est}
 \int_{B_r}\chi_\varepsilon Q^{s-2}|D^2u|^2\,dx
 \le C_{\rm cov}\int_{B_t}\chi_\varepsilon Q^s|g|^2\,dx
      +C_{\rm add}\frac{|B_R|M^s}{(t-r)^2}.
\end{equation}
The same constants are valid for the bounded-source equation when $RF_0/\varepsilon^P\le\delta_0$ for a sufficiently small, possibly nonexplicit, positive $\delta_0$ depending only on $n,p,I$.
\end{lemma}

\begin{proof}
We begin with the homogeneous case, retaining the precise comparisons needed for the bounded-source case. Set $r_{\rm cap}=(t-r)/4$. As in Lemma~\ref{lem:korn}, for $x\in B_r$ with $a_x=Q(x)>\varepsilon$ take the first upward doubling radius $l_x\in(0,r_{\rm cap}]$ capped at $r_{\rm cap}$, so $\sup_{B_{l_x}(x)}Q\le2a_x$. The choice of $\kappa_*$ and \eqref{eff:pair-holder} give $(9/10)a_x\le Q\le(11/10)a_x$ on $B_{20\kappa_*l_x}(x)$.

On this positive-gradient ball take $\psi^2u_i$ as a test function in the differentiated weak equation, where $\psi$ is the piecewise linear radial Lipschitz cutoff equal to one in $B_{5\kappa_*l_x}$ and to zero outside $B_{10\kappa_*l_x}$, with $|D\psi|\le(5\kappa_*l_x)^{-1}$. We include the bounded-source terms in this calculation. Summing over $i$ and using the upper and lower bounds for the differentiated matrix, we obtain
\begin{align*}
 \mu a_x^{p-2}\int\psi^2|D^2u|^2\,dx
 &\le 2B_*a_x^{p-2}(11/10)a_x\int\psi|D^2u||D\psi|\,dx\\
 &\quad+\sqrt n F_0\int\psi^2|D^2u|\,dx
        +2(11/10)a_xF_0\int\psi|D\psi|\,dx.
\end{align*}
The first two terms each consume at most a quarter of the left integral by Young's inequality. For the last use $2(11/10)a_xF_0\psi|D\psi|\le a_x^p|D\psi|^2+(121/100)a_x^{2-p}F_0^2\psi^2$. Since $F_0^2a_x^{4-2p}\le a_x^2/l_x^2$ under $RF_0/\varepsilon^P\le1$, we obtain
\begin{equation*}
 \int_{B_{5\kappa_*l_x}(x)}|D^2u|^2\,dx
 \le C_H a_x^2 l_x^{n-2}.
\end{equation*}
The weak testing on this positive-gradient ball is justified as in Subsection~\ref{sec:absorption}; no derivative of the source is assumed. For the homogeneous solution, simply set $F_0=0$. Hence the corresponding local truncated weighted Hessian integral is bounded by
\begin{equation}\label{eff:small-cost}
 K_{\rm loc}\chi((11/10)a_x/\varepsilon)a_x^s l_x^{n-2}.
\end{equation}

When a doubling occurs, choose $\rho_x\in[l_x,r_{\rm cap}]$ maximizing $A_x(l)^s/l^2$ with $A_x(l)=\max_{\overline B_l(x)}Q$, and set $A_x=A_x(\rho_x)$. Thus $A_x\ge2a_x$ and $a_x^s/l_x^2\le A_x^s/\rho_x^2$. If there is no doubling, or $\rho_x\ge r_{\rm cap}/L_*$, the cost in \eqref{eff:small-cost} is at most $16L_*^2K_{\rm loc}M^s l_x^n/(t-r)^2$.

In the remaining case $\rho_x<r_{\rm cap}/L_*$, the maximizing property gives $A_x(j\rho_x)\le j^{2/s_-}A_x$ for $1\le j\le L_*$. Under rescaling by the length $\rho_x$ and gradient $A_x$, the contrast hypotheses of Lemma~\ref{eff:charge} apply; in particular we have the required bound on $B_{\ell_*+1}$. Consequently, with the enlarged charge set
\begin{equation*}
 {\cal C}_x=B_{\ell_*\rho_x}(x)\cap\{Q>3A_x/4\},
\end{equation*}
we have $\int_{{\cal C}_x}Q^s|g|^2\,dx\ge(c_{\rm hom}/2)A_x^s\rho_x^{n-2}$. We have weakened the homogeneous lower bound by a factor two to use the same constant subsequently. Observe that $Q>(3/2)a_x$ in $\mathcal C_x$. The monotonicity of $\chi$ and the maximizing property therefore bound the cost \eqref{eff:small-cost} by
\begin{equation}\label{eff:charge-pay}
 \frac{2K_{\rm loc}}{c_{\rm hom}}(l_x/\rho_x)^n
                 \int_{{\cal C}_x}\chi_\varepsilon Q^s|g|^2\,dx .
\end{equation}

For $z\in{\cal C}_x$, take its own capped doubling radius $l_z$ with the same spatial cap $r_{\rm cap}$. We claim that $c_*\rho_x\le l_z\le C_*\rho_x$. Indeed $B_{\rho_x}(z)\subset B_{(\ell_*+1)\rho_x}(x)$, so $Q\le A_*A_x$ there. By the definition of $c_*$ and \eqref{eff:pair-holder}, the gradient changes from its value at $z$ by at most $(3/8)A_x$ on $B_{c_*\rho_x}(z)$. Since $Q(z)>3A_x/4$, no doubling is possible there. For the second inequality, suppose $l_z>C_*\rho_x$. Then the gradient bound on $B_{C_*\rho_x}(z)$ is $2Q(z)\le2A_*A_x$. The points $x$ and a maximum point of $Q$ in $\overline B_{\rho_x}(x)$ are both in its half-radius ball, and their mutual distance is at most $\rho_x$. The gradient difference is at most $A_x/4$, by \eqref{eff:pair-holder} and the definition of $C_*$. Their actual gradient-magnitude difference, however, is at least $A_x/2$. This contradiction proves the claim. All the balls used are inside $B_t$, since $C_*\rho_x<r_{\rm cap}$, $z\in B_{r+\ell_*r_{\rm cap}/L_*}$ and $B_{r_{\rm cap}}(z)\subset B_t$.

Select disjoint balls $B_{\kappa_*l_x}(x)$ whose fivefold dilations cover the centers under consideration. The costs of the additive class sum to the second term of \eqref{eff:cover-est} by volume disjointness. At any point $z$ belonging to several charge sets of the selected balls, the small balls of the corresponding centers are all in $B_{(\ell_*+\kappa_*)l_z/c_*}(z)$ and $\rho_x\ge l_z/C_*$, by the claim just proved. Hence
\begin{equation*}
 \sum_{x:z\in{\cal C}_x}(l_x/\rho_x)^n
 \le [C_*(\ell_*+\kappa_*)/(c_*\kappa_*)]^n=C_{\rm ov}.
\end{equation*}
Integrate \eqref{eff:charge-pay} and use this weighted overlap to obtain \eqref{eff:cover-est} in the homogeneous case.

We explain how to retain the numerical constants for the bounded-source equation, without evaluating its entire first-order modulus. The argument above uses just three comparisons from \eqref{eff:pair-holder}, at predetermined normalized scales:
\begin{itemize}
 \item on a unit ball with $\sup Q\le2$ and $Q(0)=1$, the gradient difference from the center on $\overline B_{20\kappa_*}$ is at most $1/10$;
 \item on a unit ball with $\sup Q\le A_*$ and $Q(0)\ge3/4$, the gradient difference from the center on $\overline B_{c_*}$ is at most $3/8$;
 \item on a unit ball with $\sup Q\le2A_*$, two points of $\overline B_{1/2}$ at mutual distance at most $1/C_*$ have gradient difference at most $1/4$.
\end{itemize}
For homogeneous solutions, \eqref{eff:pair-holder} gives the stronger bounds $1/30$, $3/16$ and $1/8$, respectively, since $H_0=256$. The three stated comparisons persist under sufficiently small normalized bounded sources. Otherwise there would be a sequence of normalized solutions with the prescribed gradient bound, sources tending to zero and violating points in the corresponding fixed compact interior sets. Subtract constants. The ordinary interior first-order estimates for the bounded-source equation yield a subsequence converging locally in $C^1$ to a homogeneous solution with the same bound. Pass to convergent violating points. The strict homogeneous comparisons give a contradiction in each of the three cases. This determines a positive smallness threshold without changing the fixed numerical scales.

Likewise the lower bound $c_{\rm hom}/2$ for the enlarged charge persists for small normalized bounded sources on the filled ball $B_{L_*}$ with the above record bounds on $B_{L_*}$ and $B_{\ell_*+1}$. Otherwise take a sequence of normalized solutions with sources tending to zero and exponents $s_j\in I$ tending to $s\in I$. The gradients converge locally uniformly to the gradient of a homogeneous solution, retaining the center contrast and the two bounds. On each compact subset of the open set $B_{\ell_*}\cap\{Q>3/4\}$ of the limiting solution, the weak positive-gradient Hessian estimates and the uniform convergence of the gradient coefficients give lower semicontinuity of the weighted charge. Exhaust the open set to obtain a limiting charge at most $c_{\rm hom}/2$. This contradicts Lemma~\ref{eff:charge}.

Take $\delta_0$ to be the minimum of $1$, these finitely many positive thresholds, and the common threshold fixed after Lemma~\ref{lem:korn}; retain the same notation for this smaller positive number. Each rescaling of the above covering uses gradient scale at least $\varepsilon$ and length at most $R$, including the ball centered at the charged point in the upper radius comparison. Thus $RF_0/\varepsilon^P\le\delta_0$ implies the required normalized smallness conditions. The local Hessian cost, charge lower bound, weighted overlap and additive term are unchanged, proving the lemma.
\end{proof}

\begin{theorem}[Effective strict-subthreshold exponent]\label{eff:main}
For $I=[m/2,3m/2]$ choose $s_*(n,p)$ as in \eqref{eff:main-exp}, with the constants in \eqref{eff:alpha}, \eqref{eff:lin-const}--\eqref{eff:lin-est}, \eqref{eff:warm-def}, \eqref{eff:geometry-const}--\eqref{eff:charge-const} and \eqref{eff:overlap-const}. The assertions of Theorems~\ref{thm:main} and~\ref{thm:source} hold with this explicitly specified $s_*(n,p)\in(0,m)$. In the bounded-source assertion, the source-scale multiplicative constant is not evaluated numerically.
\end{theorem}

\begin{proof}
In the production calculation of Lemma~\ref{lem:production}, keep the coefficients of the longitudinal and tangential terms. Write $A_h=\int\eta^2\chi_\varepsilon Q^sh^2$ and $A_g=\int\eta^2\chi_\varepsilon Q^s|g|^2$. The normal cutoff term is bounded by $\zeta A_h+(P^2/\zeta)\int Q^s|D\eta|^2$, and the tangential cutoff term by $(s_-/4)A_g+(4/s_-)\int Q^s|D\eta|^2$. The term $(s-2m)jh$ is bounded by another $\zeta A_h$ and a constant times the ordinary source integral. The cross term from the cutoff derivative is treated by $|jh|\le(P/2)h^2+j^2/(2P)$, consuming half of the favorable cutoff-derivative term and an additional source cost. All other source terms are lower order. This gives, directly from the identity in Lemma~\ref{lem:production},
\begin{align*}
 (3s_-/4)A_g&\le(P(m-s)+2\zeta)A_h
 +C_\zeta\int Q^s|D\eta|^2\,dx
 +C_\zeta F_0^2\varepsilon^{s-2P}|\supp\eta|.
\end{align*}
Combine with \eqref{eff:cover-est}, taking the charge-side ball to be the intermediate half-radius between the inner and outer radii. Since $Q^sh^2\le Q^{s-2}|D^2u|^2$, the leading coefficient of the resulting outer truncated Hessian integral is
\begin{equation*}
 \Theta=\frac{4C_{\rm cov}}{3s_-}\big(P(m-s)+2\zeta\big).
\end{equation*}
Take $\zeta=3s_-/(128C_{\rm cov})$ and $s=s_*$ from \eqref{eff:main-exp}. The definition of $C_{\rm cov}$ shows that $\Theta\le1/8$. The dyadic radius iteration of Subsection~\ref{sec:absorption}, which costs the factor $4^j$ at step $j$, therefore converges, since $4\Theta\le1/2$. The subsequent cutoff limit, nonlinear-gradient chain rule, local moment estimate and regularized-flux estimate at the source scale are unchanged. In the homogeneous case $\varepsilon\downarrow0$ is unrestricted. In the bounded-source case take $\varepsilon=(RF_0/\delta_0)^{1/P}$ and $\lambda=2\varepsilon$, with the common $\delta_0$ fixed in Lemma~\ref{eff:cover}. Thus $K=2\delta_0^{-1/P}$ in the effective bounded-source estimate. The exponent does not involve $\delta_0$.
\end{proof}

\begin{corollary}[Effective increment of the unweighted range]\label{eff:delta}
Let $p_c=3+2/(n-2)$, $p_\pm=p_c\pm1/2$, $P_\pm=p_\pm-1$ and $J=[p_-,p_+]$. Calculate the preceding constants as follows: take $I=[1,3]$; use $p_+$ in \eqref{eff:alpha} and in the warm coefficient bounds \eqref{eff:warm-def}; use $P_-$ in the definitions of $\nu_*,\ell_*$ in \eqref{eff:geometry-const}; and use $\mu=(9/10)^{p_+-2}$, $B_*=P_+(11/10)^{p_+-2}$ in \eqref{eff:overlap-const}. All other formulas are unchanged with these common constants. Denote the resulting covering coefficient by $C_{\rm cov,J}$ and put
\begin{equation}\label{eff:delta-formula}
 \sigma_J=\min\left\{\frac14,\frac3{64P_+C_{\rm cov,J}}\right\},
 \qquad
 \delta_n=\min\left\{\frac12,\frac{n-1}{2(n-2)}\sigma_J\right\}.
\end{equation}
Then $\delta_n>0$ depends only on $n$ and is given by finite numerical expressions. Each scalar $p$-harmonic function belongs to $W^{2,2+\varepsilon_{\rm hi}}_{\rm loc}$ for $2<p\le p_c+\delta_n$ with $\varepsilon_{\rm hi}=\varepsilon_{\rm hi}(n,p)>0$.
\end{corollary}

\begin{proof}
The common interval satisfies $1<m(p)<3<2P$ for every $p\in J$. Lemma~\ref{eff:holder-lemma} holds uniformly by the upper-endpoint ellipticity ratio. The explicit bound $M_c$ at $p_+$ and the ellipticity bounds $\lambda_{\rm ell}=q_0^{p_+-2}$, $\Lambda_{\rm ell}=P_+A^{p_+-2}$ dominate the corresponding coefficient bounds for every $p\in J$. Apply Lemma~\ref{eff:linear} directly with these common inputs in each step of the derivative bootstrap in Lemma~\ref{eff:warm}. In the Riccati argument, using $P_-$ only decreases the admissible divergence error and increases the available length. The local Hessian cost and the covering are bounded by the indicated common constants. Thus the choice $s(p)=m(p)-\sigma_J$ lies in $I$ and the absorption coefficient is at most $1/8$ uniformly over $p\in J$.

Since $m(p_c)=2$ and $m'(p)=(n-2)/(n-1)$, \eqref{eff:delta-formula} gives $s(p)\le2-\sigma_J/2$ on $[p_c,p_c+\delta_n]$. For $2<p<p_c$ the exponent in Theorem~\ref{eff:main} already satisfies $s_*<m<2$. Apply the anchored higher-integrability theorem with the corresponding $s$, and use the locally Lipschitz identity $Du=|N_s|^{2/s-1}N_s$. It follows that $Du\in W^{1,2+\varepsilon_{\rm hi}}_{\rm loc}$.
\end{proof}

We emphasize that the above argument evaluates the strict gain in the Hessian weight and the resulting positive increment $\delta_n$. The subsequent reverse H\"older improvement $\varepsilon_{\rm hi}$ in the unweighted Sobolev exponent is not evaluated here. Neither the formulas nor the argument claim the optimal weight exponent or the conjectured optimal H\"older exponent of the gradient.

\section*{Further remarks}
The proof gains an exponent by using a prescribed gradient contrast on a filled ball. The geometric lemma is essential to that gain: it could not be replaced by a pointwise inequality for arbitrary Hessian matrices satisfying the normalized trace. A related difficulty occurs in the recent work of Wang \cite[Lemma 6.3 and Section 6.2]{Wang} on monotonicity formulas for $p$-Green functions on manifolds of nonnegative Ricci curvature, where the weighted Hessian estimate with weight $|Du|^{p+\beta-4}$ is obtained for $\beta>(n-p)/(n-1)$. Putting $s=p+\beta-2$ gives $s>m$ in the Euclidean notation of this paper. The present argument instead gives a strict improvement of the local power range for scalar $p$-harmonic functions in Euclidean domains. It does not identify the optimal exponent $s_*$, nor does it establish the conjectured optimal local H\"older continuity of the gradient.


\begingroup\small
\begin{thebibliography}{99}
\bibitem{AKM}
B. Avelin, T. Kuusi and G. Mingione, \emph{Nonlinear Calder\'on--Zygmund theory in the limiting case}, Arch. Ration. Mech. Anal. \textbf{227} (2018), 663--714. \href{https://doi.org/10.1007/s00205-017-1171-7}{doi:10.1007/s00205-017-1171-7}.

\bibitem{CM}
A. Cianchi and V. Maz'ya, \emph{Second-order two-sided estimates in nonlinear elliptic problems}, Arch. Ration. Mech. Anal. \textbf{229} (2018), 569--599. \href{https://doi.org/10.1007/s00205-018-1223-7}{doi:10.1007/s00205-018-1223-7}.

\bibitem{DiB}
E. DiBenedetto, \emph{$C^{1+\alpha}$ local regularity of weak solutions of degenerate elliptic equations}, Nonlinear Anal. \textbf{7} (1983), 827--850.

\bibitem{DFZZ}
H. Dong, F. Peng, Y. R.-Y. Zhang and Y. Zhou, \emph{Hessian estimates for equations involving $p$-Laplacian via a fundamental inequality}, Adv. Math. \textbf{370} (2020), 107212. \href{https://doi.org/10.1016/j.aim.2020.107212}{doi:10.1016/j.aim.2020.107212}.

\bibitem{FMP}
M. Fogagnolo, L. Mazzieri and A. Pinamonti, \emph{Geometric aspects of $p$-capacitary potentials}, Ann. Inst. H. Poincar\'e Anal. Non Lin\'eaire \textbf{36} (2019), 1151--1179. \href{https://doi.org/10.1016/j.anihpc.2018.11.005}{doi:10.1016/j.anihpc.2018.11.005}.

\bibitem{Giaquinta}
M. Giaquinta, \emph{Multiple integrals in the calculus of variations and nonlinear elliptic systems}, Annals of Mathematics Studies, 105, Princeton University Press, Princeton, NJ, 1983.

\bibitem{HS}
A. Haarala and S. Sarsa, \emph{Global second order Sobolev-regularity of $p$-harmonic functions}, J. Differential Equations \textbf{339} (2022), 232--260. \href{https://doi.org/10.1016/j.jde.2022.08.022}{doi:10.1016/j.jde.2022.08.022}.

\bibitem{Jin}
\begingroup\raggedright
Kaiwen Jin, \emph{Second-order Sobolev regularity of three-dimensional $p$-harmonic functions}, unrefereed online manuscript, October 13, 2026. GitHub repository \texttt{yuiyuiui/AIM} (author's fork of \texttt{MColbrook/AIM}), commit \texttt{2e3c856584ed33f10b88389dd0ff65b9518e88d6}. Source file:
\url{https://github.com/yuiyuiui/AIM/blob/2e3c856584ed33f10b88389dd0ff65b9518e88d6/research/solutions/434-p-harmonic-second-sobolev-derivatives/paper/main.tex}.\par
\endgroup

\bibitem{KS}
A. Karppinen and S. Sarsa, \emph{Local second order regularity of solutions to elliptic Orlicz--Laplace equation}, Nonlinear Anal. \textbf{253} (2025), 113737. \href{https://doi.org/10.1016/j.na.2024.113737}{doi:10.1016/j.na.2024.113737}.

\bibitem{MW}
J. J. Manfredi and A. Weitsman, \emph{On the Fatou theorem for $p$-harmonic functions}, Comm. Partial Differential Equations \textbf{13} (1988), 651--668. \href{https://doi.org/10.1080/03605308808820556}{doi:10.1080/03605308808820556}.

\bibitem{MMS}
L. Montoro, L. Muglia and B. Sciunzi, \emph{Optimal second order boundary regularity for solutions to $p$-Laplace equations}, Calc. Var. Partial Differential Equations \textbf{64} (2025), 50. \href{https://doi.org/10.1007/s00526-024-02902-6}{doi:10.1007/s00526-024-02902-6}.

\bibitem{Morrey}
C. B. Morrey, Jr., \emph{Multiple Integrals in the Calculus of Variations}, Grundlehren der mathematischen Wissenschaften, 130, Springer-Verlag, Berlin, 1966.

\bibitem{Sarsa}
S. Sarsa, \emph{Note on an elementary inequality and its application to the regularity of $p$-harmonic functions}, Ann. Fenn. Math. \textbf{47} (2022), 139--153. \href{https://doi.org/10.54330/afm.112699}{doi:10.54330/afm.112699}.

\bibitem{Tolksdorf}
P. Tolksdorf, \emph{Regularity for a more general class of quasilinear elliptic equations}, J. Differential Equations \textbf{51} (1984), 126--150.
\bibitem{Wang}
Y.-Z. Wang, \emph{Elliptic monotonicity formulas for $p$-Green functions}, preprint (2026), \href{https://arxiv.org/abs/2611.32756v1}{arXiv:2611.32756v1}.

\end{thebibliography}
\endgroup
\end{document}